\documentclass[12pt]{article}
\usepackage{amsmath, amssymb}
\usepackage{mathrsfs}

\newtheorem{theorem}{Theorem}[section]
\newtheorem{lemma}{Lemma}[section]
\newtheorem{corollary}{Corollary}[section]
\newtheorem{proposition}{Proposition}[section]

\newcommand{\ad}{{\mbox{\upshape{ad}}}}
\newcommand{\adr}{{\mbox{\upshape{ad}}_r}}

\newcommand{\Bc}{{\check{B}}} 

\newcommand{\C}{{\mathbb C}}

\newcommand{\N}{{\mathbb N}}
 
\newcommand{\cC}{{\mathcal C}}
\newcommand{\Cc}{{\mathcal C}} 
\newcommand{\cF}{{\mathcal F}}
\newcommand{\cFt}{{\mathcal F}^\theta} 
\newcommand{\cG}{{\mathcal G}}
 
\newcommand{\cI}{{\mathcal I}}

\newcommand{\cM}{{\mathcal M}}

\newcommand{\qfield}{\cC} 
\newcommand{\gfrak}{{\mathfrak g}}
\newcommand{\Gr}{{\mathrm{Gr}}} 
\newcommand{\Grt}{\mathrm{Gr}^\theta}
\newcommand{\gr}{{\mathrm{gr}}} 
\newcommand{\hfrak}{{\mathfrak h}}
\newcommand{\hs}{h^\Sigma} 
\newcommand{\id}{{\mbox{id}}}

\newcommand{\kow}{{\varDelta}}

\newcommand{\alphatil}{\tilde{\alpha}}

\newcommand{\etatil}{\tilde{\eta}}
 
\newcommand{\mfrak}{{\mathfrak m}} 
\newcommand{\mutil}{\tilde{\mu}}

\newcommand{\nm}{n_{min}}

\newcommand{\ov}{{\overline{v}}}

\newcommand{\oY}{{\overline{Y}}}
 
\newcommand{\ot}{\otimes}
\newcommand{\odmu}{\overline{d_\mu}}
\newcommand{\odeta}{\overline{d_\eta}}

\newcommand{\otau}{\overline{\tau(2\mu)}}

\newcommand{\Ppi}{P(\pi)} 
\newcommand{\Ppiplus}{P^+(\pi)}

\newcommand{\rk}{\mbox{rank}}

\newcommand{\pt}{{\pi_\Theta}}
\newcommand{\Ptheta}{P_{Z(\Bc)}} 
\newcommand{\Q}{\mathbb Q}

\newcommand{\Tt}{{T_\Theta}}
 
\newcommand{\Ttc}{{\check{T}_\Theta}}

\newcommand{\uqg}{{U}}
\newcommand{\uqgdef}{{U_q(\mathfrak{g})}}

\newcommand{\Tc}{{\check{T}}} 
\newcommand{\Uc}{{\check{U}}}
\newcommand{\vep}{\varepsilon} 
\newcommand{\wght}{\mathrm{wt}}
\newcommand{\wurz}{\pi} 
\newcommand{\oYlam}{\overline{Y}_\lambda}
\newcommand{\Z}{{\mathbb Z}}

\begin{document}

\title{The center of quantum symmetric pair coideal subalgebras}

\author{Stefan Kolb\footnote{Supported by the German Research Foundation (DFG)}\, 
        and Gail Letzter\footnote{Supported by grants from the National Security Agency}\\
        \\
        {\small Mathematics Department}\\
        {\small Virginia Tech}\\
        {\small Blacksburg, VA 24061}\\
        {\small kolb@math.vt.edu,\quad letzter@math.vt.edu}}

\date{February 27, 2006}

\maketitle

\abstract{The theory of quantum symmetric pairs as developed by the
second author is based on coideal subalgebras of the quantized universal enveloping algebra
for a semisimple Lie algebra. This paper investigates the  center of these coideal 
subalgebras, proving that the center is a polynomial ring. A basis of the center
is given in terms of a submonoid of the dominant integral weights.}


\begin{center}
   {\Large Introduction}
\end{center}
After the introduction of quantum groups in the mid 1980s it was natural to look for 
quantum analogs of compact symmetric spaces. The program was motivated by a desire
to interpret I.~G.~Macdonald's large families of orthogonal polynomials 
\cite{a-Macdo00} as zonal spherical functions on such quantum analogs. 
It became evident that quantum symmetric spaces were best described in terms of one-sided 
or two-sided coideals inside the quantized universal enveloping algebra.
We refer to this characterization as the theory of quantum symmetric
pairs.  By now, there exists a well established theory of quantum
symmetric pairs and many of the initial questions have been answered.
  
There are two major approaches to the construction of quantum symmetric pairs.
In a series of papers M.~Noumi, T.~Sugitani, and M.S.~Dijkhuizen
(\cite{a-Noumi96}, \cite{a-NS95}, \cite{a-Dijk96}, \cite{a-DN98}, 
\cite{a-DS99}, and references therein) constructed them for all
classical compact symmetric spaces. The constructions were done
case-by-case and depended on solutions of the so-called reflection equation, which
were given explicitly. Taking a unified approach, the second author classified and 
analyzed quantum symmetric pairs for all compact symmetric spaces
(see \cite{MSRI-Letzter}, \cite{a-Letzter-memoirs} and references therein).
The fundamental objects, quantum symmetric pair coideal subalgebras, were described 
explicitly in terms of generators and relations \cite{a-Letzter03}. The reflection 
equation does not appear in this development.
  
By construction, quantum symmetric pair coideal subalgebras are 
quantum analogs of universal enveloping algebras of Lie subalgebras fixed by an involution
inside a semisimple Lie algebra. These fixed Lie subalgebras are reductive. One expects a 
representation theory for their quantum analogs. At present, only special
examples have been studied. For the symmetric pair of
type AI, the representations were investigated in detail in  \cite{a-IK05}. 
For the symmetric pair of type AIII, certain representations were explicitly 
constructed in \cite{a-OblSto05}. These representations were used to construct 
quantum analogs of sections of homogeneous vector bundles over Grassmann manifolds. 
This, in turn, led to new interpretations of Macdonald-Koornwinder polynomials.
  
The present paper can be viewed as the first step in developing a 
unified representation theory of quantum symmetric pair coideal subalgebras. We determine
the center of such algebras. As in the classical case, it is shown that
the center is a polynomial ring. More precisely, let $\theta$ be an 
involutive automorphism of a finite-dimensional complex semisimple Lie algebra $\gfrak$
and let $\gfrak^\theta$ denote the fixed Lie subalgebra. Let $\Bc$ 
be a corresponding quantum symmetric pair coideal subalgebra inside the simply connected 
universal enveloping algebra $\Uc$ of $\gfrak$.
  
\medskip
  
\noindent{\bf Theorem:} The center $Z(\Bc)$ is a polynomial ring in $\rk(\gfrak^\theta)$
                        variables.
  
\medskip
  
\noindent Moreover, we show that $Z(\Bc)$ has a natural basis indexed by a subset
$\Ptheta$ of the set of dominant integral weights for $\gfrak$. Just as for the ordinary
quantized enveloping algebra, we expect central characters to be 
useful in the analysis of completely reducible $\Bc$-modules and the classification of 
simple $\Bc$-modules.
  
Part of our initial motivation in studying the center $Z(\Bc)$ was to gain a better
understanding of the relation between the two main constructions of quantum symmetric pairs.
As explained in \cite[Section 6]{a-Letzter99a}, Noumi-Sugitani-Dijkhuizen type coideals are
contained in quantum symmetric pair coideal subalgebras. In the opposite direction,
it was pointed out in \cite{a-Ktobe} that any element in  $Z(\Bc)\cap (\adr \Uc)\tau(2\mu)$
leads to a solution of the reflection equation and hence associates a coideal subalgebra of
Noumi-Sugitani-Dijkhuizen type to the quantum symmetric pair coideal subalgebra $\Bc$.
Thus suitable central elements in $\Bc$ for exceptional $\gfrak$ allow
Noumi-Sugitani-Dijkhuizen type constructions
for compact symmetric spaces of exceptional type. The present paper guarantees the
existence of such central elements.
  
There are several techniques for analyzing the center of quantized universal 
enveloping algebras. M.~Rosso \cite{a-Rosso2} and, similarly,
C.~De Concini and V.G.~Kac \cite{a-DCoKa90} used an analog of the Harish-Chandra
homomorphism to describe the center. For quantum symmetric pair coideal subalgebras, 
however, it is not yet clear how to define an analog of a Cartan subalgebra. Thus a 
Harish-Chandra homomorphism is not available. In the present paper, we follow a different 
approach. A.~Joseph and the second author \cite{a-JoLet2} determined all
finite-dimensional submodules of $\Uc$ under the adjoint action. The resulting locally 
finite part $F_r(\Uc)$ contains the center of $\Uc$ in an immediately identifiable way.
It was observed in \cite{a-Letzter99b} that $Z(\Bc)$ is a subalgebra of $F_r(\Uc)$.
Hence any central element in $\Bc$ is contained in a sum of spherical submodules of
$F_r(\Uc)$. To determine $Z(\Bc)$ we exploit the structure of the direct summand
$(\adr \Uc)\tau(2\mu)$, for $\mu$ a dominant integral weight, of the locally finite
part. Moreover, using the structure of $\Bc$ and two filtrations invariant under the
adjoint action, we show that $(\adr \Uc)\tau(2\mu)$ contains an element of
$Z(\Bc)$ up to terms of lower filter degree if and only if $\mu\in \Ptheta$.
  
The situation considered here is reminiscent of an approach taken in \cite{a-Fauqu98},
\cite{a-FauquJo01} to determine the semicenter of the quantized enveloping algebra of a
parabolic subalgebra of $\gfrak$. The semicenter is also contained in the locally finite part of
$\Uc$. It is noteworthy that the condition for $(\adr \Uc)\tau(2\mu)$ to contain 
an element in the semicenter \cite[Th\'eor\`eme 3]{a-Fauqu98} resembles   the condition in
the definition (\ref{PBZ-def}) of $\Ptheta$.
  
This paper is organized in nine sections.  After fixing notation, we briefly recall the
construction and some properties of the quantum symmetric pair $\Uc$, $\Bc$.
In Section \ref{ad-inv-fil}, the filtration invariant under the adjoint action
from \cite{a-JoLet2} is compared to a similar filtration  based on the restricted
root system of the symmetric pair. It is shown that the locally finite parts of the
associated graded algebras coincide as $\Uc$-modules.
This will allow us to work with the second filtration when necessary.
  
In Section \ref{sphsub}, we obtain preliminary properties of direct summands
$(\adr \Uc)\tau(2\mu)$ of the locally finite part that contain highest
weight components of central elements in $\Bc$. A refinement of these properties leads
to the definition of $\Ptheta$ in Section \ref{HWVinGMT}.

There are two tensor product decompositions of $\Uc$, the 
triangular and the quantum Iwasawa decomposition. Using the first,
one can read off the minimal filter degree with respect to the 
filtrations introduced in Section \ref{ad-inv-fil}.
On the other hand, using the second, one easily checks whether an
element belongs to $\Bc$. In Sections \ref{triang} and \ref{qIwa}, it is shown that 
weight estimates and filter degrees can also be read off from elements decomposed with 
respect to the quantum Iwasawa decomposition.
 
The second decomposition is also at the heart of producing a projection map
from $\Uc$ onto $\Bc$ as defined in Section \ref{projUB}. Applying results from the previous
sections, we show that this projection map preserves weight estimates and filter
degrees. The projection map is then used to approximate general elements in $\Uc$
by elements in $\Bc$ with similar properties.
 
In Section \ref{HWVinGMT}, it is shown that for every $\mu\in\Ptheta$ there exists a unique
highest weight vector of spherical weight in $(\adr \Uc)\tau(2\mu)$ that also
occurs as a highest weight component of an element in $\Bc$. This highest weight vector
is the highest weight component of a central element $d_\mu\in Z(\Bc)$ constructed in
Section \ref{ZBbasis}. We prove that $\{d_\mu\,|\,\mu\in\Ptheta\}$ is a basis of $Z(\Bc)$.
   
Finally, in Section \ref{ZBpoly}, the set $\Ptheta$ is determined explicitly for every
symmetric pair. It turns out that $\Ptheta$ is a free monoid of rank $\rk(\gfrak^\theta)$.
This result is used to show that $Z(\Bc)$ is a polynomial ring in $\rk(\gfrak^\theta)$
variables.
  
For the convenience of the reader we have listed all commonly used notation -- in order of
appearance -- in an appendix.
  
\section{Preliminaries}
Let ${\C}$ denote the complex numbers, $\Z$ denote the integers, 
and $\N$ denote the nonnegative integers.   

\vspace{.5cm}
\noindent
{\bf 1.1 Complex semisimple Lie algebras.} 
  Let $\gfrak$ be a finite-dimen-sional complex semisimple Lie algebra of rank $n$ and 
  $\hfrak$ be a fixed Cartan subalgebra of $\gfrak$. Let $\Delta$ denote the root system
  associated with $(\gfrak,\hfrak)$. Choose an ordered basis 
  $\wurz=\{\alpha_1,\dots,\alpha_n\}$ of simple roots for $\Delta$.
  Identify $\hfrak$ with its dual via the Killing form. The induced
  nondegenerate symmetric bilinear form on $\hfrak^*$
  is denoted by $(\cdot,\cdot)$. Let $W$ denote the Weyl group associated to 
  the root system $\Delta$ and let $w_0$ denote the longest element
  in $W$ with respect to $\wurz$.
  We write $Q(\wurz)$ for the root lattice and $P(\pi)$ for the 
  weight lattice
  associated to the root system $\Delta$.
  Set $Q^+(\pi)=\N\pi$ and let $P^+(\pi)$ be the set of dominant integral weights with 
  respect to $\pi$. We will denote the fundamental weights in $P^+(\pi)$
  by $\omega_1.\dots,\omega_n$.
  Let $\leq$ denote the standard partial ordering on 
  $\hfrak^*$.  In particular, $\mu\leq\gamma$ if and only  if $\gamma-\mu\in Q^+(\pi)$.

\vspace{.5cm}
 \noindent
 {\bf 1.2 The quantized enveloping algebra.}
 Let $q$ be an indeterminate and let $M$ be a positive integer so 
 that $(M/2)(\lambda,\mu)\in \Z$ for all $\mu\in P(\pi)$. Let $q^{1/M}$ be 
 a fixed $M^{th}$ root of $q$ in the algebraic closure of $\C(q)$.
 Set ${\cal C}$ equal to the field of rational functions $\C(q^{1/M})$. 
 We consider here the $q$-deformed universal enveloping algebra $\uqgdef$ as the 
 $\qfield$-algebra generated by elements $t_i,t_i^{-1}$, $x_i$, $y_i$, 
 where $i=1,\dots,n$, subject to the  relations  as given in \cite[Section 3.2.9]{b-Joseph}.
 In particular, one has 
 \begin{enumerate}
   \item[(i)] The  $t_1,\dots, t_n$ generate a free abelian group of 
   rank $n$.
   \item[(ii)]$t_ix_j=q^{(\alpha_i,\alpha_j)}x_j t_i$ and $
       t_iy_j=q^{-(\alpha_i,\alpha_j)}y_j t_i$   for   all $1\leq 
       i\leq n$.
   \item[(iii)] $x_iy_j-y_jx_i=\delta_{ij} 
		      (t_i-t_i^{-1})/(q^{(\alpha_i,\alpha_i)/2}-q^{-(\alpha_i,\alpha_i)/2})$ for   all 
		      $1\leq i,j\leq n$.
   \item[(iv)]The $ x_1,\dots, x_n$ and the $y_1,\dots, y_n$ satisfy
              the quantum Serre relations (see for example \cite[(1.7)]{MSRI-Letzter}).
 \end{enumerate}   
  We often write $U=\uqgdef$.
  Note that it is necessary to work over the field $\qfield$ 
  instead of $\C(q)$ so that $U$ may be enlarged to the simply 
  connected quantized enveloping algebra and the extension of the 
  quantized enveloping algebra used in Section 5.1.  
       
  The algebra $U$ also has a Hopf algebra structure with coproduct 
  $\kow$, antipode $\sigma$, and counit $\epsilon$ satisfying
  \begin{align*}
    \kow x_i&=x_i\ot 1 +t_i\ot x_i, & \sigma (x_i)&=-t_i^{-1}x_i,& \epsilon(x_i)&=0,\\
    \kow y_i&=y_i\ot t_i^{-1}+ 1\ot y_i, & \sigma (y_i)&=-y_it_i,&\epsilon(y_i)&=0,\\
    \kow t_i&=t_i\ot t_{i},& \sigma(t_i)&=t_i^{-1},& \epsilon(t_i)&=1.\\ 
  \end{align*}
  The algebra $U$ acts on itself from the right by the right adjoint
  action, denoted by $\adr$.  This action is completely determined by 
  the action of the generators of $U$ as follows:
  \begin{align}
    (\adr t_i)b&=t_i^{-1}bt_i,\nonumber\\ 
    (\adr x_i)b&=t_i^{-1}b x_i-t_i^{-1}x_ib\label{ad},\\ 
    (\adr y_i)b&=b y_i-y_it_ib t_i^{-1}\nonumber 
  \end{align}
  for all $b\in U$ and $1\leq i\leq n$. Throughout this paper 
  we will encounter various right $U$-modules for which the right 
  action is induced by the adjoint action of $U$ on itself. We 
  will frequently call such modules $(\adr U)$-modules.   
  Similarly, we will talk of $(\adr M)$-modules for subalgebras $M$ 
  of $U$ and its extensions.
       
  Let $U^+$ denote the subalgebra  of $U$ generated by  $x_1,\dots, x_n$ and let $U^-$ 
  denote the subalgebra of $U$ generated by $y_1,\dots,y_n$.
  Similarly, let $G^+$ and $G^-$ denote the subalgebras
  of $U$ generated by $x_1t_1^{-1},\dots, x_n t_n^{-1}$ and 
  $y_1t_1,\dots,y_nt_n$, respectively. Let $T$  denote the
  multiplicative group generated by $t_1^{\pm 1},\dots, t_n^{\pm 1}$ and let 
  $U^0$ denote the subalgebra of $U$ generated by $T$.
  The algebra $U$ admits a triangular decomposition. In particular, 
  the multiplication map induces an isomorphism
  \begin{align}\label{trig-decomp}
    U^-\ot U^0\ot U^+\rightarrow U.
  \end{align}
  The same holds if one or both of $U^-$ and $U^+$
  are replaced by $G^-$ and $G^+$, respectively.
  
Let $\tau$ denote the group isomorphism from $Q(\pi)$ onto $T$ given 
by $\tau(\alpha_i)=t_i$ for $i= 1,\dots, n$.  We can enlarge $T$ to a 
group $\check T$ such that $\tau$ extends to an isomorphism from 
$P(\pi)$ onto $\check T$.  
It is often convenient to consider the simply connected quantized 
enveloping algebra $\Uc$ generated by $U$ and $\check T$.   
In particular, $\Uc$ is  generated by 
$x_i,y_i$, for $i=1,\dots, n$, and $\tau(\lambda),$ for $\lambda\in 
P(\pi)$, satisfying all the relations of $U$ and $\check T$ as well as 
\begin{align}\label{taurelns}
  \tau(\lambda)x_j=q^{(\lambda,\alpha_j)}x_j \tau(\lambda),\qquad
  \tau(\lambda)y_j=q^{-(\lambda,\alpha_j)}y_j \tau(\lambda)
\end{align}
for all $\lambda\in P(\pi)$ and $1\leq j\leq n$. 
The subalgebra of $\Uc$ generated by the elements $\tau(\lambda)$,
$\lambda\in \Ppi$, is denoted by $\Uc^0$.

The following notation is used throughout this paper.
For any multiindex $I=(i_1,\dots,i_m)$, $1\le i_j\le n$, define
$(yt)_I=y_{i_1}t_{i_1}\cdots y_{i_m}t_{i_m}$.
Set ${\rm wt}(I)=\alpha_{i_1}+\dots+\alpha_{i_m}$. Note that $(yt)_I\in G^-$.
Let $\cI$ be a set of multiindices such that $\{(yt)_I\,|\,I\in \cI\}$
is a basis of $G^-$.  
   
\vspace{.5cm}
\noindent{\bf 1.3 Left and right $U$-modules.}
Let $M$ be a left $U$-module.   Given a vector subspace $S$ of $M$
and a weight $\mu\in P(\pi)$, the $\mu$ weight space of $S$ is the
subspace $S_{\mu}$ defined by
\begin{align}\label{wspaceone} 
  S_{\mu}=\{a\in S\,|\,t_i\,a =q^{(\alpha_i,\mu)}a\ \mbox{for all } i=1,\dots,n\}.
\end{align}
A weight vector $v$ of $S$ is called a highest weight vector
provided that $v\neq 0$ and $x_i v=0$ for all $1\leq i\leq n$.   Similarly, a
lowest weight vector $v$ is a nonzero weight vector with $y_i v=0$ for
$1\leq i\leq n$.  
    
Following standard notation as in \cite{b-Joseph}, the weight spaces
of $\Uc$ are  weight spaces with respect to the left
adjoint action.  In particular,   given a subset $M$ of $\Uc$
and a weight $\mu\in P(\pi)$, the  subset $M_{\mu}$ is defined by
\begin{align}\label{wspacetwo} 
  M_{\mu}=\{a\in M\,|\,\tau(\lambda)a\tau(\lambda)^{-1}=q^{(\lambda,\mu)}a\ 
    \mbox{for all }\lambda\in P(\pi)\}.
\end{align}  
For $\mu\in \Ppiplus$, let $V(\mu)$ denote the uniquely determined
finite-dimen-sional simple left $\uqg$-module of highest weight $\mu$.
    
To stick to the conventions of \cite{MSRI-Letzter} we will consider left coideal subalgebras
of $ \Uc$ in Section 1.5. It follows from Theorem \ref{locfinBthm} (below) that it will be natural to
consider right $U$-modules in order to determine the center of these coideal 
subalgebras. Therefore
we will consider the dual space $V(\mu)^*$ always with its natural right $U$-module 
structure defined by $fu(v)=f(uv)$ for all $f\in V(\mu)^*$, $v\in V(\mu)$, and $u\in  U$.
    
For a subspace of a right $U$-module $M$ one defines the right $\mu$ weight space
as in (\ref{wspaceone}) with $t_ia$ replaced by $at_i$.  
Note that the right $U$-module  $V(\mu)^\ast$ has a nonzero $\mu$ weight space. 
However, the elements of this weight space are annihilated by each of the $y_i$,
$1\leq i\leq n$. Hence we call a nonzero weight vector $v$ of a right $U$-module a highest
weight vector if $vy_i=0$ for all $i=1,\dots,n$. 
It should also be noted that $M_\mu$ as defined in (\ref{wspacetwo}) is the $-\mu$
weight space for the right adjoint action. 
   
   \vspace{.5cm}
   \noindent
   {\bf 1.4 Classical infinitesimal symmetric pairs.} 
   Let $\theta$ be an 
   involutive Lie algebra automorphism of $\gfrak$. Write 
   $\gfrak^\theta$ for the Lie subalgebra of $\gfrak$ consisting of 
   elements fixed by $\theta$. The pair $(\gfrak,\gfrak^\theta)$ is 
   a (infinitesimal) 
   symmetric pair. A classification of involutions and symmetric pairs up 
   to isomorphism can be found in 
   \cite[Chapter X, Sections 2, 5, and 6]{b-Helga78} and 
   \cite[Section 4.1.4]{b-OV94}.  We further assume that the 
   involution $\theta$ is maximally split with respect to $\hfrak$ and the 
   chosen triangular decomposition of $\gfrak$ in the sense of 
   \cite[Section 7]{MSRI-Letzter}. It should be noted that every 
   involution of $\gfrak$ is conjugate to a maximally split one inside the 
   group of Lie algebra automorphisms of $\gfrak$.
   
   The involution $\theta$ induces an involution $\Theta$ of the root 
   system $\Delta$ of $\gfrak$ and thus an automorphism of $\hfrak^\ast$. 
   Define $\pt=\{\alpha_i\in \wurz\,|\,\Theta(\alpha_i)=\alpha_i\}$.  
   Let $p$ be the permutation of $\{1,\dots,n\}$ such that
   \begin{align}\label{pdefn}
     \Theta(\alpha_i)\in -\alpha_{p(i)}+\Z\pi_{\Theta}\ \mbox{for all}\ \alpha_i\notin \pt
   \end{align}  
   and $p(i)=i$ if $\alpha_i\in \pt$. Moreover, let $\pi^\ast$ be the set 
   of all $\alpha_i\in \pi\setminus\pt$ such that $i=p(i)$ or $i<p(i)$. 

   There is a second root system associated to the symmetric pair 
   $(\gfrak, \gfrak^{\theta})$ defined as follows.  Given $\alpha\in \hfrak^\ast$, set 
   $\alphatil=(\alpha-\Theta(\alpha))/2$. The subset 
   \begin{align*}
     \Sigma=\{\alphatil\,|\,\alpha\in\Delta \,\mbox{and}\, 
			       \Theta(\alpha)\neq \alpha\}
   \end{align*}
   of $\hfrak^\ast$ is the restricted root system associated to the pair 
   $(\gfrak,\gfrak^\theta)$. The set of simple roots of $\Sigma$ is just
   $\{\alphatil_i\,|\,\alpha_i\in \pi^\ast\}$. Let $P(\Sigma)$ denote the 
   weight lattice of $\Sigma$. Let $P^+(\Sigma)$ be the set of dominant integral  
   weights with respect to $\Sigma$. 
              
   The partial ordering $\leq$ on $\hfrak^*$ restricts to a partial ordering on 
   $P(\Sigma)$.  It is also useful to introduce a restricted version of the standard partial 
   ordering.  In particular, given $\mu,\gamma\in P(\Sigma)$, we write
   $\mu\leq_r\gamma$ if $\gamma-\mu\in \sum_{\alpha_i\in \pi^*}\N\tilde\alpha_i$.

\vspace{.5cm}
\noindent{\bf 1.5 Quantum symmetric pairs.} The main result of 
\cite{MSRI-Letzter} is the classification for maximally split symmetric pairs 
$(\gfrak,\gfrak^\theta)$ of all left coideal subalgebras $B$ inside $\uqg$ satisfying the 
following two properties:
\begin{align*}
  \bullet\,&\mbox{$B$ specializes to $U(\gfrak^\theta)$ as $q$ goes to $1$.}\\
  \bullet\,&\mbox{\parbox[t]{11cm}{If $B\subseteq C$ and $C$ is a 
      subalgebra of $\uqg$ which specializes to $U(\gfrak^\theta)$ then $B=C$.}}
\end{align*}
The pair $U$, $B$ is called a quantum symmetric pair and is an analog of the pair of
enveloping algebras $U(\gfrak)$, $U(\gfrak^\theta)$. We briefly review the structure
of the coideal subalgebras $B$. 
 
Let $\cM$  denote the Hopf subalgebra of $\uqg$ generated by 
$x_i$, $y_i$, $t_i^{\pm 1}$ for $\alpha_i\in \pt$. Define the subgroup $\Tt$ of $T$ by
\begin{align*}
  \Tt=\{\tau(\lambda)\,|\,\lambda \in Q(\pi)\,\mbox{and}\,           
			\Theta(\lambda)=\lambda \}
\end{align*}
The subalgebra $B$ of $ \uqg$ is generated by $\cM$, $\Tt$, and 
elements $B_i$ for $\alpha_i\in \pi\setminus\pt$. For the explicit form 
of the elements $B_i$ consult \cite[Section 7]{MSRI-Letzter}. In this paper we only use that
\begin{align}\label{Bi-def}
  B_i=y_i t_i + d_i\tilde{\theta}(y_i)t_i+s_i t_i
\end{align}
where $d_i,s_i\in \qfield$, $\tilde{\theta}(y_i)\in G^+_{\Theta(-\alpha_i)}$, and
$\tilde{\theta}(y_i)t_i\in \cM x_{p(i)}\Tt$.  
It is shown in \cite{MSRI-Letzter} that  any maximal left coideal
subalgebra  of $\uqg$ that specializes to $U(\gfrak^{\theta})$ 
as $q$ goes to $1$  is isomorphic to an algebra $B$ of the above form via 
a Hopf algebra automorphism of $\uqg$. 
   
Set $B_i=y_it_i$ for $\alpha_i\in \pi_{\Theta}$ and set $\cM^+=\cM\cap U^+$. 
Recall the definition of  $\cI$ given in Section 1.2. Given a multiindex 
$I=(i_1,\dots, i_m)$ in $\cI$, set $B_I=B_{i_1}\cdots B_{i_m}$.
By \cite[Section 7]{MSRI-Letzter}, we have the direct sum decomposition	    
\begin{align}\label{Bdirectsum}
  B=\oplus_{I\in \cI}B_I\cM^+\Tt.
\end{align}
It follows that when $b\in B$ is written as a sum of weight vectors in 
$\check U$, any lowest weight term is an element of $G^-{\cal M}^+T_{\Theta}$.
	
Let $T'_{\Theta}$ be the subalgebra of $\check U$ defined by 
\begin{align*} 
  T'_{\Theta}=\{\tau(\lambda)|\ \lambda\in P(\pi){\rm \ and \ }
	 \Theta(\lambda)=\lambda\}.
\end{align*}
Let $\check B$ be the subalgebra of $\Uc$ generated by $B$ and $T'_{\Theta}$.  
It follows from properties of $B$ that $\check B$ is a left 
coideal subalgebra of $\check U$ that specializes to $U(\gfrak^{\theta})$ as $q$ goes to $1$.
It should be noted that	the direct sum decomposition (\ref{Bdirectsum}) as 
well as the subsequent remark also hold for $\check B$ if $T_{\Theta}$ is 
replaced by the group $T'_{\Theta}$.

\vspace{.5cm}
\noindent
{\bf 1.6 The locally finite part of ${ \check U}$.}
  The locally finite part  of $\check U$ with respect to  the right 
  adjoint action of $ U$ on $\check U$ is the set 
  \begin{align*}
    F_r(\check U)= \{a\in \check U|\ \dim((\adr U)a)<\infty\}.
  \end{align*}  
  We recall here some basic facts about $F_r(\check U)$ from \cite{a-JoLet2}
  (see also \cite{b-Joseph}). However, it should be noted that \cite{a-JoLet2}
  is written in terms of the left adjoint action while here we  focus on the right adjoint 
  action. Definitions and theorems from \cite{a-JoLet2} will be translated accordingly. 
				 
  The vector space $F_r(\check U)$ is a subalgebra of $\check U$.
  As an $(\adr U)$-module, it can be decomposed as follows \cite[Thm.~4.10]{a-JoLet2}
  \begin{align}\label{ldecomp}
    F_r(\Uc)=\bigoplus_{\mu\in\Ppiplus}(\adr U)\tau(2\mu).
  \end{align}
  Recall further the isomorphism $(\adr U)\tau(2\mu)\cong V(\mu)^\ast\ot V(-w_0\mu)^\ast$ 
  of right $U$-modules \cite[Corollary 3.5]{a-JoLet2}.
			
  Consider $\mu\in P^+(\pi)$. The following  formula is useful in analyzing the 
  $(\adr U)$-module generated by $\tau(2\mu)$.
  Let $\zeta$ and $\xi$ be elements in $Q^+(\pi)$.  The description of the 
  right adjoint action combined with the relations of $U$ ensure that 
  \begin{align}\label{adUzetaxi}
    (\adr U^-_{-\zeta}U^+_{\xi})\tau(2\mu)\in \sum_{\alpha\leq \zeta}\, 
    \sum_{\beta\leq \xi}\,\ \sum_{\makebox[1cm]
    {$0\le\eta\le \zeta-\alpha \atop 0\le \eta\le\xi-\beta$}}
      U^-_{-\alpha}G^+_{\beta}\tau(2\mu-2\eta).
  \end{align}

Let $Z(\Bc)$ denote the center of $\Bc$.  The fact that $\Bc$ 
is a left coideal implies that any element in $Z(\Bc)$ spans
a one-dimensional $(\adr B)$-invariant submodule of $\Uc$.   
The following observation will allow us to use the well known 
structure of the locally finite part $F_r(\check U)$ in the investigation of the center 
of $\Bc$.
\setcounter{theorem}{5}
\begin{theorem}\label{locfinBthm}  \cite[Theorem 1.2]{a-Letzter-memoirs} 
  Each finite-dimensional $(\adr B)$-submodule  of $\check U$ 
  is a subset  of $F_r(\check U)$. In particular, $Z(\Bc)$ is a 
  subalgebra of $F_r(\check U)$.
\end{theorem}
		       
\vspace{.5cm}
\noindent{\bf 1.7 Spherical modules.}
Given a $U$-module $M$, let $M^B$ denote the set of $B$-invariant 
elements inside of $M$.  For example, if $M$ is a left $U$-module, 
then $M^B=\{v\in M|\ bv=\vep(b)v\ \mbox{for all }b\in M\}$.   
A finite-dimensional simple   $\uqg$-module $V(\lambda)$, $\lambda\in P^+(\pi)$, is 
called spherical if $V(\lambda)^B$ is one dimensional.
It is shown in \cite[Theorem 3.4]{a-Letzter03}
that  $V(\lambda)$ is a finite-dimensional spherical $\uqg$-module  if and only if 
$\lambda\in 2P^+(\Sigma)$. If $\lambda\in 2P^+(\Sigma)$, then we refer to $\lambda$ as a 
spherical weight. 
		      
Fix $\lambda\in P^+(\pi)$ and let $v_{\lambda}$ denote the highest weight  vector of 
$V(\lambda)$.  Assume further that $V(\lambda)$ is spherical and let $v^B$ be a nonzero 
$B$-invariant element of $V(\lambda)$.  Rescaling if necessary, we can write 
\begin{align}\label{vBweights}
  v^B= v_{\lambda}+\sum_{\{\mu\in P(\Sigma)|\ \mu<_r\lambda\}}w_{\mu}
\end{align}
where each $w_{\mu}$ is a vector of weight $\mu$ inside 
$V(\lambda)$ \cite[Theorem 3.6]{a-Letzter03}.  

By \cite[Theorem 7.6]{MSRI-Letzter} the left $B$-module $V(\lambda)$ is 
a direct sum of irreducible left $B$-modules. Hence the right
$B$-module $V(\lambda)^\ast$ also decomposes into irreducible right
$B$-modules. Thus $V(\lambda)^\ast$ contains a right
$B$-invariant element if and only if $V(\lambda)$ contains a left
$B$-invariant element. Moreover, this holds if and only if $\lambda$
is spherical.
		     
\section{Two $(\adr U)$-invariant filtrations}\label{ad-inv-fil}
  
In this section, we consider two $(\adr U)$-invariant filtrations on 
$\check U$.   The first filtration  ${\cal F}$ is  one of the basic tools used to analyze the locally finite part of 
$\check U$ in  \cite{a-JoLet2} (see also 
\cite[Chapter 7]{b-Joseph}). The second filtration ${\cal F}^{\theta}$ is a 
modification of ${\cal F}$ that is based on the restricted root 
system $\Sigma$ in the same way as ${\cal F}$ is based on the root 
system $\pi$.  An isomorphism of the locally finite parts of the 
graded algebras  of $\check U$ with respect to these two different filtrations is 
established at the end of the section.

\vspace{.5cm}	 
\noindent
{\bf 2.1 The standard $(\adr U)$-invariant filtration.}
Let $h:\Q\pi\rightarrow \Q$ denote the height function defined by
\begin{align}\label{h-def}
  h(\sum_i n_i\alpha_i)=\sum_i n_i.
\end{align} Note that  the image of $P(\pi)$ under the height function  is contained in $\frac{1}{N}\Z$ for 
some sufficiently large positive integer $N$.  Consider the ${{1}\over{N}}\Z$-filtration 
$\cF$ of $\check U$ defined by the following degree function
\begin{align}\label{filtdef}
  \deg(x_it_i^{-1})=\deg(y_i)=0,\quad \deg \tau(\lambda)=h(\lambda) 
\end{align}
for all $i=1,\dots,n$ and  $\lambda\in P(\pi)$.
Note that ${\cal F}$ restricts to a $\Z$-filtration on $\uqg$.
Define $\Gr:=\gr_\cF(\Uc)$, $\Gr^0:=\gr_\cF(\Uc^0)$, 
$\Gr^-:=\gr_\cF(U^-)$, and $\Gr^+:=\gr_\cF(G^+)$ to be the associated graded algebras.
The description (\ref{ad}) of 
the right adjoint action of $U$ on $\check U$  implies that 
$\cF$ is an $(\adr U)$-invariant filtration. In particular, $\Gr$ 
inherits the structure of an $(\adr U)$-module from $\Uc$. Given
$a\in \Uc$, we denote its graded image in $\Gr$ by $\bar a$.  More 
precisely, for each  $a\in {\cal F}_n(\check U)\setminus {\cal F}_{n-1/N}(\check U)$ we 
write $\bar a$ for the image of $a$ in the subspace   
${\cal F}_n(\check U)/{\cal F}_{n-1/N}(\check U)$ of $\Gr$.
  
Note that $U^-\subseteq \cF_0U^-$ and $\cF_{-1/N}U^-=\{0\}$.
Hence $U^-\cong \Gr^-\cong\Gr_0^-$. Similarly one obtains $G^+\cong \Gr^+\cong \Gr^+_0$. 
Note  that $\Gr^0$, $\Gr^-$, and $\Gr^+$ are subalgebras of $\Gr$. By the triangular 
decomposition (\ref{trig-decomp}) the multiplication map 
\begin{align*}
  \Gr^-\ot\Gr^0\ot\Gr^+\rightarrow \Gr
\end{align*}
is an isomorphism. Using this isomorphism, we may identify 
$\Gr^-\ot\Gr^0$ with $\Gr^-\Gr^0$ and $\Gr^0\ot\Gr^+$ with 
$\Gr^0\Gr^+$ as vector spaces.   
  
\vspace{.5cm}
\noindent
{\bf 2.2 The locally finite part of $\Gr$.} 
Given $\mu\in P^+(\pi)$, define the subspace  $K(2\mu)^+$ of $\Gr^+$  and the subspace 
$K(2\mu)^-$ of $\Gr^-$  by 
\begin{align}
  (\adr U^-)\overline{\tau(2\mu)}&=K(2\mu)^-\ot \qfield\overline{\tau(2\mu)},
                               \label{Klambdadef}\\
  (\adr U^+)\overline{\tau(2\mu)}&=\qfield\overline{\tau(2\mu)}\ot K(2\mu)^+.
                               \label{Klambdadef2}
\end{align}  
By  \cite[Section 4.9]{a-JoLet2} one has    
\begin{align}
  (\adr U)\overline{\tau(2\mu)}=K(2\mu)^-\ot  \qfield\overline{\tau(2\mu)} \ot K(2\mu)^+.
    \label{grinto}
\end{align}

\setcounter{lemma}{1}

\begin{lemma}\label{adra} 
  Given $\mu\in P^+(\pi)$ and $G\in K(2\mu)^-_{\alpha}\ot 
  \qfield\overline{\tau(2\mu)}\ot K(2\mu)^+_{-\beta}$, there 
  exists an element $a\in\sum_{\zeta\leq \alpha}\sum_{\xi\leq \beta} U^-_{-\zeta}U^+_{\xi}$
  such that $(\adr a){\overline{\tau(2\mu)}}=G$.
\end{lemma}
  
\noindent
  {\bf Proof:} Without loss of generality, we may assume that 
  $G=u'\ot \overline{\tau(2\mu)}\ot v'$ for some $u'\in K(2\mu)^-_{\alpha}$ and 
  $v'\in K(2\mu)^+_{-\beta}$. By (\ref{Klambdadef}) and (\ref{Klambdadef2}) there 
  exist $u\in U^-$ and $v\in U^+$ such that $(\adr u)\tau(2\mu)=u'\ot\overline{\tau(2\mu)}$ 
  and $(\adr v)=\overline{\tau(2\mu)}\ot v'$.  We may further assume that $u\in 
  U^-_{-\alpha}$ and $v\in U^+_{\beta}$ because the adjoint action on 
  $\overline{\tau(2\mu)}$ preserves weight spaces.    
  Using the fact that $K(2\mu)^-$ is an $(\adr U^+)$-submodule of $\Gr^-$  
  \cite[Section 4.9]{a-JoLet2} combined with the relations of $U$  one obtains 
  \begin{align*}
    (\adr uv)&\otau=(\adr v)(u'\ot \otau)\\
     &\in q^{(\alpha,\beta)}u'\ot \otau\ot v' +
     \sum_{\makebox[0cm]{$\eta< \alpha \atop \xi<\beta$}}
     K(2\mu)^-_{\eta}\ot \otau \ot K(2\mu)^+_{-\xi}.
  \end{align*}
  The lemma follows by induction on the height of $\alpha$ and $\beta$. 
  $\blacksquare$

\vspace{.5cm}

By \cite[Lemma 7.1.4]{b-Joseph} the $(\adr U)$-modules $(\adr U)\overline{\tau(2\mu)}$ 
and $(\adr U)\tau(2\mu)$ are isomorphic for each $\mu\in P^+(\pi)$.   Moreover, the locally 
finite part $F_r(\Gr)$ of $\Gr$ is isomorphic to $F_r(\check U)$ as 
an $(\adr U)$-module \cite[Sec. ~4.10]{a-JoLet2}, \cite[Sec. 7.1]{b-Joseph}.   
This isomorphism is given by the function $g:F_r(\Uc)\rightarrow F_r(\Gr)$ defined by 
\begin{align}
  g((\adr u)\tau(2\mu))=(\adr u)\overline{\tau(2\mu)}\label{griso}
\end{align} 
for all $\mu\in \Ppiplus$ and $u\in U$.
The isomorphism $g$ combined with the decomposition of $F_r(\check U)$ in (\ref{ldecomp}) 
yields the following direct sum decomposition:
\begin{align}\label{grlf}
  F_r(\Gr)=\bigoplus_{\mu\in P^+(\pi)}(\adr U)\overline{\tau(2\mu)}.
\end{align}
  
\vspace{.5cm} \noindent
  {\bf 2.3 A coarse $\adr$-invariant filtration.}
  The second filtration is defined using a height function 
  associated to the restricted root system instead of the ordinary 
  height function.   More precisely, the height function $h^{\Sigma}$ 
  associated to the restricted root system is the function 
  $h^{\Sigma}:\Q\Sigma\rightarrow\Q$ defined by  
  \begin{align}\label{restrictedheight}
    h^{\Sigma}(\sum_{\alpha_i\in \pi^\ast} n_i\alphatil_i)=\sum_{\alpha_i\in \pi^\ast}n_i.
  \end{align}
  The function $h^{\Sigma}$ can be  extended to a function on $\Q\pi$ by setting 
  $h^{\Sigma}(\lambda)=h^{\Sigma}(\tilde\lambda)$ for all $\lambda\in \Q\pi$. 
  Note that $h^{\Sigma}$ restricts to a function from $P(\pi)$ to ${{1}\over{N}}\Z$ for a 
  sufficiently large positive integer $N$.  For the remainder of this paper, we 
  assume that $N$ has been chosen so that both height functions $h$ 
  and $h^{\Sigma}$ map $P(\pi)$ to ${{1}\over{N}}\Z$.
	   
  Define a ${{1}\over{N}}\Z$-filtration $\cFt$ of    $\check U$ by the degree function
  \begin{align}\label{filtthetadefn}
    \deg(x_it_i^{-1})=\deg(y_i)=0,\qquad
    \deg(\tau(\lambda))=\hs(\lambda)
  \end{align}
  for all $i=1,\dots, n$ and $\lambda\in P(\pi)$. Note that $\cFt$ 
  restricts to a $\Z$-filtration on $U$.  Moreover,   $\cFt$ is 
  an $(\adr U)$-invariant on both $U$ and $\check U$. 
  The analogous filtration invariant under the left adjoint
  action has been used in \cite[Sec.~5]{a-Letzter04}.
  Let $\Grt=\gr_{{\cal F}^{\theta}}(\Uc)$ be the associated graded algebra. Given
  $a\in \Uc$, we denote its graded image in $\Grt$ by  $\bar {a}^{\theta}$.  
     
  Consider weights $\alpha$ and $\beta$ in $Q(\pi)$ such that $\alpha\leq \beta$.
  Note that if $h^{\Sigma}(\alpha)<h^{\Sigma}(\beta)$ then $\alpha<\beta$ and 
  so $h(\alpha)<h(\beta)$. This observation will be useful in comparing the two filtrations 
  ${\cal F}$ and ${\cal F}^{\theta}$.
     
  \vspace{.5cm}
  \noindent
  {\bf 2.4 The locally finite part of $\Grt$.} 
  The locally finite part of $\Grt$ with respect to the right adjoint action of $U$ on $\check U$ 
  is the set
  \begin{align*}
    F_r(\Grt)=\{a\in \Grt\,|\,\dim((\adr U)a)<\infty \}.
  \end{align*}
  An immediate consequence of the following result is that the $\frac{1}{N}\Z$-graded
  $(\adr U)$-module $F_r(\Grt)$ can be obtained from the $\frac{1}{N}\Z$-graded 
  $(\adr U)$-module $F_r(\Gr)$ by renumbering
  degrees of the submodules of the form $(\adr U)\overline{\tau(2\mu))}$.
     
  \setcounter{proposition}{3}   
  \begin{proposition}\label{FGrtheta}
    The map $g^\theta:F_r(\Uc)\rightarrow F_r(\Grt)$ defined by 
    \begin{align*}
      g^{\theta}((\adr u)\tau(2\mu))=(\adr u)\overline{\tau(2\mu)}^\theta
    \end{align*}
    for all $\mu\in P^+(\pi)$ and $u\in U$ is an isomorphism of right $U$-modules.
  \end{proposition}
  {\bf Proof:} 
    Recall that $F_r(\Uc)$ and $F_r(\Gr)$ are isomorphic as 
    $(\adr U)$-modules via the map $g$ defined in (\ref{griso}). 
    Note that by the discussion  at the end of Section 2.3,
    one gets $\ker g^\theta\subseteq \ker g$. Therefore $g^\theta$ is injective.
     
    To prove surjectivity let $v\in \cFt_s\Uc$ represent an element
    $\ov^\theta\in F_r(\Grt_s)$.  In particular, 
    $\dim((\adr U)\ov^\theta)<\infty$.
    Without loss of generality, we may assume that there exist $\lambda_1,\dots, \lambda_k$ such that 
    $h^{\Sigma}(\lambda_i)=s$ for all $i=1,\dots, k$ and 
    \begin{align*}
      v\in \sum_{i=1}^k U^-  \tau(\lambda_i)G^+.
    \end{align*}
     It follows  from the description of the right adjoint action (\ref{ad}) and the 
     defining relations of $U$ that 
    \begin{align*} 
      (\adr U)  v\in \sum_{i=1}^k\sum_{\mu\leq\lambda_i}U^-\tau(\mu)G^+.
    \end{align*}
    Since $(\adr U)\bar v^{\theta}$ is finite dimensional, there exist subspaces $U_i$ of 
    $U^-$ and $G_i$ of $G^+$ such that
    \begin{align*}
      (\adr U)  v\subseteq \sum_{i=1}^k  U_i
      \tau(\lambda_i)G_i+\sum_{i=1}^k
      \sum_{\mu<\lambda_i}U^-\tau(\mu)G^+.
    \end{align*}
		  
    Let $m$ be the smallest 
    rational number such that $v\in \cF_m(\check U)$.   
    It follows that 
    \begin{align}\label{v-form}
      (\adr U)v\subseteq \sum_{\{i|h(\lambda_i)=m\}} U_i
       \tau(\lambda_i)G_i+\sum_{\{\mu\in \Ppi| h(\mu)<m\}}U^-\tau(\mu)G^+.
    \end{align}
    Note that (\ref{v-form}) 
    ensures that  $(\adr U)\bar v$  is contained in the image of the  
    finite-dimensional subspace 
    $\sum_{\{i|  h(\lambda_i)=m\}}U_i\tau(\lambda_i)G_i$ inside $\Gr$. In particular, 
    $\bar v\in F_r(\Gr)$.   It follows  from (\ref{grlf}) that 
    \begin{align}\label{vfour-form}
      v\in\sum_{\{i|  h(\lambda_i)=m\}}(\adr U)\tau(\lambda_i)
      +\sum_{\{\mu\in \Ppi|  h(\mu)<m\}}U^-\tau(\mu)G^+.
    \end{align}
    Moreover,  $\lambda_i\in 2P^+(\pi)$ for all $i$ such that $h(\lambda_i)=m$. 
    In particular, $m$ is a nonnegative element of ${{1}\over{N}}\Z$.
	  
    Using (\ref{vfour-form}), we can choose $w\in 
    \sum_{h(\lambda_i)=m}(\adr U)\tau(\lambda_i)$
    such that $v-w$ is of smaller filter degree than $m$ with respect
    to  ${\cal F}$.  Since each $\lambda_i\in 2P^+(\pi)$, we see that $\bar w^{\theta}$ 
    is in $F_r(\Grt)$. Assume first that $m=0$. Suppose that 
    $\bar v^{\theta}\neq \bar w^{\theta}$. It follows that $\bar v^{\theta}-\bar 
    w^{\theta}=\overline{(v-w)}^{\theta}$. Hence, $\overline{(v-w)}^{\theta}$ is an element 
    of $F_r(\Grt)$. Arguing as above with $v$ replaced by $v-w$, we see that 
    $v-w\in {\cal F}_c(\check U)$ for some  nonnegative element $c$ of ${{1}\over{N}}\Z$. 
    Moreover, $c$ must be strictly less than $m$ and $m=0$.
    This contradiction forces $\bar v^{\theta}=\bar w^{\theta}$.  
    Therefore $\bar v^{\theta}=g^{\theta}(w)$ and so $\bar v^{\theta}$ is in 
    the image of $g^{\theta}$. Now assume that $m>0$.
    If $\bar v^{\theta}=\bar w^{\theta}$, then once again $\bar v^{\theta}$ is in 
    the image of $g^{\theta}$.  Thus we may assume that
    $\bar v^{\theta}\neq\bar w^{\theta}$.  Replacing $v$ by 
    $v-w$, we obtain $v-w\in {\cal F}_c(\check U)$ for 
    some  nonnegative element $c$ of ${{1}\over{N}}\Z$ strictly less than $m$.
    By induction, $\overline{v-w}^{\theta}\in g^{\theta}(F_r(\check U))$.   
    The lemma now follows from the fact that $\overline{w}^{\theta}$ 
    is in the image of $g^{\theta}$ and 
    $\bar v^{\theta}-\bar w^{\theta}=\overline{v-w}^{\theta}$.
  $\blacksquare$
      
  \vspace{.5cm}
  \noindent
  {\bf 2.5 $B$-invariant elements of $\Gr$ and $\Gr^{\theta}$.}  
  Recall that $\check U$ is semisimple as an $(\adr T)$-module. Hence both
  $\Gr$ and $\Gr^{\theta}$ are semisimple $(\adr T)$-modules. 
  It follows from \cite[Lemma 3.4]{a-Letzter99b} that any element $b\in \Gr^B$ 
  generates a finite-dimensional $(\adr U)$-submodule. Thus $\Gr^B\subseteq F_r(\Gr).$
  Similarly one has $(\Gr^{\theta})^B\subseteq F_r(\Grt).$
  The second part of the following result is an immediate consequence of 
  Theorem \ref{locfinBthm}, Proposition \ref{FGrtheta}, and the 
  fact that the map $g$ defined by (\ref{ldecomp}) is an isomorphism. 
  \setcounter{theorem}{4}
   \begin{theorem}\label{GrB} 
    We have the following inclusions of $B$-modules:
    \begin{align*}
      \Gr^B\subseteq F_r(\Gr),\qquad(\Gr^{\theta})^B\subseteq F_r(\Grt).
    \end{align*}
    Moreover, the map $g$ restricts to an isomorphism of $\check U^B$ onto $\Gr^B$ and
    the map $g^{\theta}$ restricts to an isomorphism of $\check U^B$
    onto $(\Gr^{\theta})^B$.
  \end{theorem}
  
\section{Spherical submodules of $(\adr U)\tau(2\mu)$}\label{sphsub}
  The first lemma of this section describes a 
  property satisfied by  highest weight vectors inside submodules of the 
  form $(\adr U)\tau(2\mu)$.  We then   
  establish a multiplicity result concerning  certain spherical  finite-dimensional 
  simple $U$-modules 
  inside $(\adr U)\tau(2\mu)$, for $\mu\in P^+(\pi)$.
  This latter result   will help us   determine  which 
  submodules $(\adr U)\tau(2\mu)$ of $F_r(\Uc)$ contain 
  nonzero elements of $Z(\Bc)$ up to lower degree terms with respect 
  to the filtration ${\cal F}^{\theta}$. 

\vspace{.5cm}
   
\noindent
{\bf 3.1 Highest weight vectors of $\Gr$.} 
  By  \cite[4.8]{a-JoLet2}, $(\adr U^-){\overline{\tau(2\mu)}}$ contains a unique 
  one-dimensional highest weight space with respect to the right 
  adjoint action.   Moreover, the weight of this space is 
  $\mu-w_0\mu$ and all other weights of $(\adr U^-)\overline{{\tau(2\mu)}}$ are 
  strictly smaller than $\mu-w_0\mu$. The next lemma provides a useful 
  property of highest weight vectors in $\Gr$.
  \setcounter{lemma}{0}
  \begin{lemma}\label{hwcor}
    Let $\mu\in P^+(\pi)$ and   let $v$ be a highest weight vector in 
    $(\adr U)\overline{\tau(2\mu)}$. Then there exist nonzero weight vectors 
    $X\in K(2\mu)^+$ and $Y\in K(2\mu)^-_{\mu-w_0\mu}$
    such that 
    \begin{align*} 
      v\in Y\otimes \overline{\tau(2\mu)}\otimes X+\sum_{\eta<\mu-w_0 \mu} 
      K(2\mu)^-_\eta\ot\qfield\overline{\tau(2\mu)}\ot K(2\mu)^+.
    \end{align*}
  \end{lemma}
		 
  \noindent
  {\bf Proof:} We can write $v=\sum_{i=1}^kY_i\otimes X_i$ where 
  each $Y_i$ is a weight vector  of $K(2\mu)^-\ot\overline{\tau(2\mu)}$
  and each $X_i$ is a weight vector  of $K(2\mu)^+$ with respect
  to  the right adjoint action.  Assume further that   $\{X_i|\ 1\leq i\leq k\}$ is a 
  linearly independent set. Let $\lambda$ denote the weight of $v$ and let $\eta_i$ 
  denote the weight of $Y_i$. Choose $s\in\{1,\dots, k\}$ such that $\eta_s$ is maximal. 
  Then for all $j=1,\dots,n$ one has
  \begin{align*}
    0&=(\adr y_j)v\\
     &=\sum_{i=1}^k (\adr y_j)Y_i\ot (\adr t_i^{-1})X_i+\sum_{i=1}^k Y_i \ot (\adr y_i)X_i\\
     &\in \sum_{\eta_i=\eta_s} (\adr y_j)Y_i\ot (\adr t_i^{-1})X_i+
       \sum_{\eta\neq \eta_s}K(2\mu)^-\ot\qfield\overline{\tau(2\mu)}\ot 
       K(2\mu)^+_{\lambda-\eta}					       
  \end{align*}
  and hence $(\adr y_j)Y_i=0$ for all $j=1,\dots,n$ and for all $i$ such that 
  $\eta_i=\eta_s$. In particular, $Y_i$ is a highest weight vector in 
  $(\adr U^-)\overline{\tau(2\mu)}$ whenever
  $\eta_i=\eta_s$.  The result now follows from the discussion preceding 
  the lemma.
  $\blacksquare$
	
\vspace{.5cm}
\noindent {\bf 3.2 Multiplicity of spherical modules.} 
Recall that $\pi_{\Theta}$ is the subset of $\pi$ consisting of those 
simple roots $\alpha_i$ such that $\Theta(\alpha_i)=\alpha_i$.  
Write $W(\pt)$ for the subgroup of the Weyl group $W$ generated by the simple reflections
corresponding to the simple roots in $\pt$. Note that $W(\pt)$ is the Weyl 
group  corresponding to the root system generated by $\pt$. Let $w_0'$ denote the
longest element in $W(\pt)$. 

Let $\mfrak$ denote the semisimple Lie subalgebra of $\gfrak 
$ generated by the positive and negative root vectors corresponding to 
the simple roots in $\pi_{\Theta}$.  Note that ${\cal M}$ can be 
identified with the quantized enveloping algebra of $\mfrak$.  Suppose that $\mu\in 
P^+(\pi)$.   Then $\mu$ restricts to a dominant integral weight with 
respect to the root system generated by $\pt$.   It follows as in the 
discussion preceding Lemma \ref{hwcor}, that the subspace of $(\adr \cM^+)\otau$ 
annihilated by each $(\adr x_i)$ for $\alpha_i\in \pt$ is 
one dimensional.   Moreover, this space is a weight space of weight 
$-\mu+w_0'\mu$ with respect to the right adjoint action of $U$.
\setcounter{theorem}{1}
\begin{theorem}\label{sphersubprop}
Let $\lambda$ be a spherical weight and let $\mu\in \Ppiplus$ such that
$\mu-w_0\mu-\lambda\in \N \pt$. Then
\begin{align}\label{mult-le}
  [(\adr U)\tau(2\mu):V(\lambda)^\ast]\le 1.
\end{align}
Moreover, 
\begin{align}\label{mult-e}
  [(\adr U)\tau(2\mu):V(\lambda)^\ast]= 1
\end{align}
if and only if $\lambda=w_0'\mu-w_0\mu$.
\end{theorem}

  \noindent{\bf Proof:} 
  As noted in Section 2.2, the right $U$-modules $(\adr U)\otau$ and $(\adr U)\tau(2\mu)$ 
  are isomorphic.  Thus it is 
  sufficient to prove (\ref{mult-le}) and (\ref{mult-e})
  using $(\adr U)\otau$ instead of $(\adr U)\tau(2\mu)$.   Let  $Y_{\lambda}$ denote a 
  highest weight vector of weight $\lambda$ inside $(\adr U)\otau$.
  By the discussion preceding Lemma \ref{hwcor}, we can choose a unique
  (up to scalar multiplication) nonzero $Y\in K(2\mu)^-$ such that $Y\otimes\otau$
  is a highest weight vector of weight   $\mu-w_0\mu$ in $(\adr U^-)\otau$. 
  Lemma \ref{hwcor} ensures that there exists a nonzero  $X\in K(2\mu)^+$ 
  such that 
  \begin{align*}
     Y_\lambda\in Y\ot \otau\ot X +\sum_{0\le\eta<\mu-w_0\mu}
     K(2\mu)^-_\eta\ot \otau\ot K(2\mu)^+_{\lambda-\eta}.
  \end{align*}
  We claim that
  \begin{align}\label{MHW}
    (\adr x_i)(\otau\ot X)=0  \mbox{ for all } \alpha_i\in \pi_\Theta.
  \end{align}
  Indeed, recall that $\lambda$ is spherical and $Y_\lambda$ is a highest weight 
  vector of weight $\lambda$. It follows that $Y_{\lambda}$ is 
  $(\adr {\cal M}T_{\Theta})$-invariant.   Hence, for $\alpha_i\in \pi_\Theta$ 
  one has
  \begin{align*}
    0=(\adr x_i)Y_\lambda\in\ ((\adr t_i) Y)&\ot((\adr x_i)(\otau\ot X))\\
      &+\sum_{0\le\eta<\mu-w_0\mu} K(2\mu)^-_\eta\ot \otau\ot K(2\mu)^+
  \end{align*}
  which proves (\ref{MHW}). 
     
  Note that $Y\otimes \otau\otimes X$ is a vector of weight $\lambda$.
  It follows that 
  $$\otau\ot X\in \otau\ot K(2\mu)^+_{-\beta} $$
  where $\beta= \mu-w_0\mu -\lambda$.  By assumption, $\beta\in {\N}\pi_{\Theta}$.
  Hence $\otau\ot X\in (\adr \cM^+)\otau$. By the discussion preceding 
  the lemma, $\otau\ot X$ spans the unique  one-dimensional subspace  
  in $(\adr \cM^+)\otau$  annihilated by all 
  $(\adr x_i)$ for $\alpha_i\in \pi_\Theta$.
  Hence up to multiplication by a nonzero scalar, any highest weight vector 
  $Y'_{\lambda}$ of weight $\lambda$ in $(\adr U)\otau$ satisfies
  \begin{align*}
    Y'_{\lambda}\in Y\ot \otau\ot X +\sum_{0\le\eta<\mu-w_0\mu} K(2\mu)^-_\eta\ot \otau\ot 
    K(2\mu)^+.
  \end{align*}
  Therefore any two highest weight vectors of weight $\lambda$ in
  $(\adr U)\otau$ are linearly dependent. Hence (\ref{mult-le}) is proved.

  By (\ref{MHW}) and the discussion preceding the lemma, the weight $-\beta$ of  $\otau\ot X
  $ satisfies $\beta=\mu-w_0'\mu$. On the other hand, 
  $\beta$ is defined to be the weight $\mu-w_0\mu-\lambda$.  Thus
  \begin{align}\label{ww'cond}
    \lambda=w_0'\mu-w_0\mu
   \end{align}
  is a necessary condition for (\ref{mult-e}) to hold.
   
  Assume now that $\lambda=w_0'\mu-w_0\mu$ holds.  Recall the isomorphism of right 
  $U$-modules $(\adr U)\tau(2\mu)\cong V(\mu)^*\ot V(-w_0\mu)^*$ noted in Section 1.6.
  To complete the 
  proof, it suffices to show that the left $U$-module $V(\mu)\ot V(-w_0\mu)$ contains a 
  highest weight vector of weight $\lambda=w_0'\mu-w_0\mu$.
  
  Let $M$ denote the ${\cal M}$-submodule of $ V(\mu)$ generated by 
  the highest weight vector $v_{\mu}\in V(\mu)$.  Similarly, let 
  $M^\ast$ denote the ${\cal M}$-submodule of $ V(-w_0\mu)$ generated by the 
  highest weight vector $v_{-w_0\mu}\in V(-w_0\mu)$. As $\lambda=w_0'\mu-w_0\mu$
  and $(\lambda,\alpha_i)=0$
  for all $\alpha_i\in \pt$ the ${\cal M}$-modules $M$ and $M^\ast$ are
  dual to each other. Hence the subset $M\ot M^\ast$ of $ V(\mu)\ot V(-w_0\mu)$
  contains a nonzero $\cM$-invariant element $v$ of weight
  $w_0'\mu-w_0\mu=\lambda$. But as $\lambda+\alpha_i\nleqslant \mu-w_0\mu$
  for all $\alpha_i\notin\pt$, the element $v$ is a highest weight vector for the action 
  of $U$ on $V(\mu)\ot V(-w_0\mu)$. This completes the proof that the condition
  $\lambda=w_0'\mu-w_0\mu$ implies (\ref{mult-e}).
   $\blacksquare$

 \noindent	   
 \section{The triangular decomposition} \label{triang}
 
The aim of this section is to extract information about elements in $B$ and $F_r(\check U)$
with respect to  the triangular decomposition $\check U\cong G^-\ot U^+\ot \check U^0$.
These results will be  used  to express certain  elements in 
$F_r(\Uc)$ as elements of $\check B$ up to terms of  lower $\cFt$ filter degree.    
  
\vspace{.5cm}
    
\noindent{\bf 4.1 Triangular decomposition of elements in $\check B$.}
Recall the partial ordering $\geq_r$ defined in Section 1.4. Set
\begin{align*}
  T_{\ge}={\rm span}\{\tau(\delta)|\, \delta\in P(\pi){\rm \  and\ }
  \tilde\delta\geq_r 0\}.
\end{align*}
Note that $T_{\geq}$ is a subalgebra of $\check U^0$ which contains 
${\cal C}[t_i| \alpha_i\in \pi]$ and $T'_\Theta$.
\setcounter{lemma}{0}
\begin{lemma}\label{BtrigLemone}
  The algebra   $\check B$ is a subset of $G^-U^+T_{\ge}$.
\end{lemma}
    
\noindent{\bf Proof:}
  The relations of $U$ guarantee that the
  subalgebra of $U$ generated by $G^-$, $U^+$, $T_\Theta'$, and
  ${\cal C}[t_i|\alpha_i\in \pi]$ is a subset of $G^-U^+T_\ge$.
  Note that by definition (\ref{Bi-def}) each $B_i$ lies in the subalgebra
  generated by $G^-$, $U^+$, $\Tt$ and ${\cal C}[t_i|\alpha_i\in \pi]$.
  Moreover, ${\cal M}\subseteq G^-U^+\Tt$. This completes the proof of the lemma.
  $\blacksquare$
   
\vspace{.5cm}
\noindent{\bf 4.2 Triangular decomposition of elements in $(\adr U)\tau(2\mu)$.}
  It follows from (\ref{ldecomp})  that ${\cal F}_m(F_r(\check U))$ is finite dimensional for
  all $m\geq 0$. Unfortunately, a similar statement does not hold for 
  the filtration ${\cal F}^{\theta}$.   In order to prove the main results of this paper, 
  it is necessary to restrict to finite-dimensional subspaces of 
  ${\cal F}^{\theta}_m(F_r(\check U))$.  This is done by intersecting 
  with spaces $G({\le}\mu)$ defined as follows. For $\mu\in P(\pi)$ set
  \begin{align}\label{Glessthanmu}
    G({\le}\mu):= U^-G^+\tau(\mu)\qfield[t_1^{-2},\dots,t_n^{-2}]
	        =\sum_{\eta \in Q^+(\pi)}U^-G^+\tau(\mu-2\eta).
  \end{align}
  Note that the sets $G({\le}\mu)$ for $\mu\in P(\pi)$ satisfy the following multiplicative 
  property:  
  \begin{align}\label{mult}
    G({\le}\mu_1)G({\le}\mu_2)=G({\le}(\mu_1+\mu_2)){\rm \ for\ all\ } 
    \mu_1,\mu_2\in P(\pi).
  \end{align}
\setcounter{proposition}{1}	 
\begin{proposition}\label{BtrigLem}
  Fix  $\mu\in P^+(\pi)$. Suppose that $v\in (ad_r U)\tau(2\mu)$ and $w\in G^-U^+ T_{\geq}$ 
  satisfy
  \begin{align}\label{v-b-inFtm-1}
    v-w \in \sum_{\alpha\not\geq 2\mu}G({\le}\alpha). 
  \end{align} 
  Then $v\in G^-U^+T_\ge$.
\end{proposition}
	 
\noindent
  {\bf Proof:} 
  The description of the right adjoint action (\ref{ad}) implies that $(\adr U)\tau(2\mu)$
  is a subset of $G({\le}2\mu)$. In particular, there exists $X\in U^-G^+\tau(2\mu)$ 
  such that
  \begin{align}\label{v-one}
    v\,&\in\,X+\sum_{\eta> 0}U^-G^+\tau(2\mu-2\eta).
  \end{align} 
  It now follows from (\ref{v-b-inFtm-1}) that 
  \begin{align*}
    w\in X+\sum_{ \alpha\neq 2\mu} U^-G^+\tau(\alpha).
  \end{align*}
  On the other hand, since $w\in G^-U^+ T_{\geq}$ we have
  \begin{align*}
    w\in \sum_{\xi,\zeta}\sum_{\tilde\gamma\ge_r 0} G^-_{-\zeta} U^+_\xi\tau(\gamma)=
    \sum_{\zeta,\xi }\sum_{\tilde\gamma-\tilde\zeta-\tilde\xi\geq_r 0}U^-_{-\zeta} G^+_\xi
    \tau(\gamma)    
  \end{align*}
  where $\zeta$ and $\xi$ run over elements in $Q^+(\pi)$ and each
  $\gamma$ is in $P(\pi)$.  It follows that 
  \begin{align}\label{Xformula}
    X\in  \sum_{\makebox[0cm]{$\zeta, \xi \atop 2\tilde\mu\geq_r \tilde\zeta+\tilde\xi$}}
    U^-_{-\zeta} G^+_\xi\tau(2\mu).
  \end{align}
By (\ref{Xformula}) and Lemma \ref{adra}, there exists 
\begin{align*}
  a\in \sum_{\makebox[0cm]{$\zeta, \xi \atop 2\tilde\mu\geq_r \tilde\zeta+\tilde\xi$}}
  U^-_{-\zeta}U^+_{\xi}
\end{align*}
such that $v=(\adr a)\tau(2\mu)$. Now (\ref{adUzetaxi}) implies that 
\begin{align*}
  v&\in \sum_{\alpha,\beta}\, \sum_{2\tilde\eta\leq_r 
  2\tilde\mu-\tilde\alpha-\tilde\beta}U^-_{-\alpha}G^+_{\beta}\tau(2\mu-2\eta)\\
  &=\sum_{\alpha,\beta}\,  \sum_{ 2\tilde\eta\leq_r 
  2\tilde\mu-\tilde\alpha-\tilde\beta}G^-_{-\alpha}U^+_{\beta}\tau(2\mu-2\eta
  -\alpha-\beta)\\
  &\subseteq \sum_{\alpha,\beta}\sum_{\{\delta\in P(\pi)|\ \tilde\delta\geq_r 0\}}
  G^-_{-\alpha}U^+_{\beta}\tau(\delta).
\end{align*}
The proposition follows from the definition of $T_{\geq}$.
$\blacksquare$
	
\vspace{.5cm}
	
\noindent
{\bf 4.3 Triangular decomposition of elements in ${\cal F}_m^{\theta}(\check U)$.}
Suppose that  $a\in {\cal F}^{\theta}_m(F_r(\check U))$ for some $m\in \frac{1}{N}\Z$.
It follows from Proposition \ref{FGrtheta} that 
\begin{align*}
  a\in \sum_{\{\gamma\in P^+(\pi)|\ h^{\Sigma}(2\gamma)\leq m\}}(\adr U)\tau(2\gamma).
\end{align*} 
Now suppose that $\mu\in P^+(\pi)$ such that $h^{\Sigma}(2\mu)>m$.
It further follows that $a\in \sum_{\alpha\not\geq 2\mu}G({\le}\alpha)$.
The next result can be viewed as an extension of
Proposition \ref{BtrigLem} with $(\adr U)\tau(2\mu)$ replaced by  $F_r(\check U)$.
\setcounter{corollary}{2}	
\begin{corollary}\label{v-bcor}
  Suppose that $w\in G^-U^+T_{\geq}$ and 
  \begin{align*}
    v\in \sum_{\{\gamma\in P^+(\pi)|\ h^{\Sigma}(2\gamma)>m\}}(\adr U)\tau(2\gamma)
  \end{align*} 
  such that $v-w\in {\cal F}^{\theta}_m( \check U)$ for some $m\in \frac{1}{N}\Z$. 
  Then $v\in G^-U^+ T_{\geq}$.
\end{corollary}

\noindent
{\bf Proof:}
We can find distinct weights $\mu_i\in P^+(\pi)$, for $i=1,\dots, r$, and 
elements $v_i\in (\adr U)\tau(2\mu_i)$ such that $v=v_1+\dots +v_r$.
Furthermore, we may assume  that $h^{\Sigma}(2\mu_i)\geq m+{{1}\over{N}}$ for each $1\leq 
i\leq r$. Reorder the $\mu_i$ if necessary so that $\mu_i>\mu_j$ implies $i>j$.
It follows that $v_i\in \sum_{\alpha\not\geq 2\mu_r}G({\le}\alpha)$ for 
$i=1,\dots, r{-}1$. Moreover, ${\cal F}_m^{\theta}(\check U)$ is a subset of 
$\sum_{\alpha\not\geq 2\mu_r}G({\le}\alpha)$ and hence 
$v_r-w\in\sum_{\alpha\not\geq 2\mu_r}G({\le}\alpha)$. An application of the preceding 
proposition yields $v_r\in G^-U^+ T_{\geq}$. The corollary now 
follows by induction on $r$ after replacing $w$ by $w-v_r$ and $v$ by $v-v_r$.
$\blacksquare$

\section{The  quantum Iwasawa decomposition}\label{qIwa}

There is a second tensor product decomposition of $\check U$ 
involving the subalgebra $\check B$ called the quantum Iwasawa 
decomposition.  Since we are interested 
in elements of $\check B$, it is often better to express elements 
in terms of the quantum Iwasawa decomposition.   This section provides 
ways to translate information between the Iwasawa and the triangular 
decompositions. 

\vspace{.5cm}
\noindent
{\bf 5.1 The  quantum Iwasawa  decomposition of $\check U$.} 
  In general, the set $\widetilde{\Ppi}=\{\tilde{\alpha}\,|\,\alpha\in {P(\pi)}\}$ 
  is not a subset of $P(\pi)$.
  As in Section 1.2, we enlarge $\Tc$ to a group $\Ttc\check{\cal A}$ such that the 
  isomorphism
  $\tau$ extends to an isomorphism from $\Ppi+\widetilde{\Ppi}$ to 
  $\Ttc\check{\cal A}$. Set $\check{\cal A}=\{\tau(\tilde\alpha)|\ \alpha\in P(\pi)\}$
  and $\Ttc=\{\tau((\alpha+\Theta(\alpha))/2)|\alpha\in P(\pi)\}$.
  Note that $\Ttc\cap \check  T=T'_{\Theta}$ and that
  $\Ttc \check{\cal A}$ is isomorphic to the group
  $\Ttc \times\check{\cal A}$. As in Section 1.2, we can form the 
  algebra $\check U\Ttc$ by insisting that relation 
  (\ref{taurelns}) holds for all $\tau(\lambda)\in \Ttc$. 
  
  Let ${\rm ad}$ denote the left adjoint action \cite[(1.2) and Section 5]{MSRI-Letzter}.
  Set $N^+$ equal to the subalgebra of $U^+$ generated by the $(\ad \cM^+)$-module
  $(\ad \cM^+)(\qfield[x_i\,|\,\alpha_i\notin \pt])$.
  The following isomorphism is a form of the quantum Iwasawa decomposition and easily 
  deduced from \cite[ Section 4]{a-Letzter-memoirs}. The multiplication map
  induces an isomorphism
  \begin{align}\label{Iwasawa}
    B\Ttc \ot \qfield[\check{\cal A}]\ot N^+\cong \check U\Ttc. 
  \end{align}
  This implies in particular that 
  $\check U=(BN^+\Ttc\check{\cal A})\cap \check U=BN^+\check U^0$.  
 
  The next lemma provides more detailed information concerning
  the triangular decomposition of elements $B_I$, for $I\in {\cal I}$.
  \setcounter{lemma}{0}
  \begin{lemma}\label{BIdecomp}  
    For all $I\in \cI$, we have 
    \begin{align*} 
      B_I\in (yt)_I+\sum_{\makebox[0cm]{$\eta<{\rm wt}(I)\atop \hs(\eta)<\hs(\wght(I))$}} 
        \,G^-_\eta U^+ T.
    \end{align*}  
  \end{lemma}
	
 \noindent
 {\bf Proof:}  This is a straightforward induction on $h({\rm wt}(I))$. 
 $\blacksquare$

 \vspace{.5cm}
 \noindent
 {\bf 5.2  Comparing the two tensor product decompositions.} 
 Estimates on the filter degree or weight of an element $u\in \Uc$ are usually
 obtained by decomposing $u$ with respect to the triangular
 decomposition (\ref{trig-decomp}). A consequence of the following lemma
 is that the summands of the
 decomposition of $u$ with respect to the quantum Iwasawa decomposition
 often satisfy the same estimates.
 \setcounter{proposition}{1}
 \begin{proposition}\label{triangIwa} 
   Assume $\lambda\in P(\pi)$ and $k\in {{1}\over{N}}\Z$. 
   Then the following relations hold.
   \begin{enumerate}
     \item[(i)]	$G({\le}\lambda)=\sum_{I\in {\cal I}}B_I[G^+\check U^0\cap
                 G({\le}(\lambda-{\rm wt}(I)))]$.
     \item[(ii)] ${\cal F}_k^{\theta}(\check U)=\sum_{I\in {\cal I}}
             B_I[G^+\check U^0\cap {\cal F}_{k-h^{\Sigma}({\rm wt}(I))}^{\theta}(\check U)]$.
     \item[(iii)] $\sum_{h^{\Sigma}(\eta)<k}\check U_{-\eta}=\sum_{I\in {\cal I}}
         B_I[G^+\check U^0\cap \sum_{h^{\Sigma}(\eta+{\rm wt}(I)))<k}\check U_{-\eta}].$
   \end{enumerate}
 \end{proposition}
   
 \noindent {\bf Proof:} 
   Note that for all three statements it suffices to show that the left hand side is  
   contained in the right hand side. We first prove (i).
   Let $a\in G({\le}\lambda)$. Without loss of generality we can assume that $a=(yt)_I a_I$
   for some $I\in \cI$ and $a_I\in U^+\Tc$. We have to show that $(yt)_Ia_I$ is also 
   contained in the set on the right hand side of (i). Note that $a\in G({\le}\lambda)$ 
   implies $a_I\in G({\le}(\lambda-\wght(I)))$. The multiplicative property 
   (\ref{mult}) implies $B_Ia_I\in G({\le}\lambda)$. By Lemma \ref{BIdecomp} one gets
   \begin{align*}
     (yt)_I a_I-B_Ia_I\in \sum_{\eta<\wght(I)} G^-_{-\eta}U^+T\cap G({\le}\lambda).
   \end{align*}
   Claim (i) now follows by induction on the height of $\wght(I)$.
   
   One obtains assertions (ii) and (iii) using a similar argument with $G({\le}\lambda)$
   replaced by $\cF^\theta_k(\Uc)$ and $\sum_{h^{\Sigma}(\eta)<k}\check U_{-\eta}$,
   respectively.
   $\blacksquare$
	 
 \vspace{.5cm}
	    
\noindent 
\section{A projection of $\check U$ onto $\check B$}\label{projUB}
 
This section defines a projection $b$ of $\check U$ onto $\check B$ 
using  the quantum Iwasawa decomposition. We further show that $b(u)$ inherits estimates
of filter degree and weight from $u$. The map $b$ is a critical tool in the construction of 
a particularly nice basis of the center of $\check B$ in Section \ref{ZBbasis}. 
   
\vspace{.5cm}
\noindent
{\bf 6.1 An $({\adr} {\cal M}T_{\Theta})$-submodule of $\check U\Ttc$.}
Let $N^+_+$ denote the intersection of $N^+$ with the augmentation ideal of $U$.
Set $\check{\cal A}^{+}:=\check{\cal A}\setminus\{1\}$. 
Let $\cG$ denote the subset   $\check{\cal A}N^{+}_++\check{\cal A}^+$ of $\check U\Ttc$.  
\setcounter{lemma}{0}
\begin{lemma}\label{IwaLemma}
  The vector space $B\Ttc  \cG$ is an  $(\adr \cM\Tt)$-module. Moreover,
  $B\Ttc \cG\cap  B\Ttc =\{0\}$.
\end{lemma}
{\bf Proof:} The second statement follows immediately from the 
  quantum Iwasawa decomposition (\ref{Iwasawa}).
      
  To prove the first statement note first that $B\Ttc \cG$ is 
  $(\adr \Tt)$-invariant. Moreover, $(\adr m)t'\in\cM t'$ for all
  $m\in \cM$, $t'\in \check{\cal A}^+$ and $\check{\cal A}^+\cM=\cM 
  \check{\cal A}^+$.
  Thus it remains to prove that 
  \begin{align*}
    (\adr x_i) N^{+}_+\subseteq \cM N^{+}_{+},\qquad (\adr y_i) N^{+}_+\subseteq \cM N^{+}_{+}
  \end{align*}
  for all $\alpha_i\in \pt$.
  In view of (\ref{ad}) these relations are equivalent to
  \begin{align*}
    N^{+}_+x_i\subset \cM N^{+}_+,\qquad N^{+}_+y_i\subset \cM N^{+}_+
  \end{align*}
  for all $\alpha_i\in \pt$.
  The first formula follows immediately from the fact that the multiplication map
  gives an isomorphism $U^+\cong \cM^+\ot N^+$ of vector spaces \cite{a-Kebe99}. 
  For the second formula one has to commute $y_i$ to the left side using the defining 
  relations (iii) of $U$ before applying $U^+\cong \cM^+\ot N^+$.
  $\blacksquare$ 
  
        
\vspace{.5cm}
    
\noindent
{\bf 6.2 A  ${\cal M}T_{\Theta}$-module decomposition.}   
Recall that $ B\Ttc $ contains the Hopf algebra ${\cal M}T_{\Theta}$.  It follows that 
$B\Ttc $ is an $(\adr {\cal M}T_{\Theta})$-module. 
Thus the quantum Iwasawa decomposition (\ref{Iwasawa})  and Lemma \ref{IwaLemma} 
yield the following direct sum 
decomposition 
\begin{align}
  \label{Iwoplus} \check U\Ttc =  B\Ttc {\cal G}\oplus B\Ttc  
\end{align} 
of $\check U\Ttc $ into ${\cal M}T_{\Theta}$-modules with respect to the right adjoint 
action. Given $u\in  \check U\Ttc $, let $b(u)$ be the projection of $u$ 
onto $\check B\Ttc$ with respect to (\ref{Iwoplus}). 
\setcounter{lemma}{1}
\begin{lemma}
  The projection $b$ restricts to a mapping of $\check U$ 
  onto $\check B$.
\end{lemma}

\noindent
{\bf Proof:}  The quantum Iwasawa decomposition for $\Uc$ below (\ref{Iwasawa}) 
implies that $b(\check U)=b(B\check U^0)$. Any element $u\in 
B\check U^0$ can be uniquely decomposed in the form $u=u_1+u_2$ where 
\begin{align*}
  u_1\in \sum_{\Theta(\lambda)\neq \lambda}B\tau(\lambda), {\qquad}u_2\in 
  \sum_{\Theta(\lambda)=\lambda}B\tau(\lambda).
\end{align*}
The claim of the lemma follows because $b(u)=u_2$ and $u_2\in \check B$. 
$\blacksquare$

\vspace{.5cm}
\noindent
{\bf 6.3 Properties of the map $b$.}
Consider $u\in \check U$. The next result shows that  
$b(u)$ retains certain characteristics of $u$ related to   degree and weight. 
\setcounter{proposition}{2}
\begin{proposition}\label{constructingLem2} 
  Assume $\lambda\in P(\pi)$ and $k\in {{1}\over{N}}\Z$. The map $b$ has the following 
  properties.
  \begin{enumerate}
    \item[(i)] $b(G({\le}\lambda))\subseteq G({\le}\lambda)$.
    \item[(ii)] $b({\cal F}_k^{\theta}(\check U))\subseteq {\cal F}_k^{\theta}(\check U)$. 
    \item[(iii)] $b(\sum_{h^{\Sigma}(\eta)<k}\check U_{-\eta})\subseteq 
        \sum_{h^{\Sigma}(\eta)<k}\check U_{-\eta}$.
  \end{enumerate}
\end{proposition}
    
\noindent{\bf Proof:} We prove assertion (i) using Proposition 
  \ref{triangIwa}(i). The other two assertions follow in a similar manner from
  Proposition \ref{triangIwa}(ii) and (iii). Assume $u\in G({\le}\lambda)$. We can write
  \begin{align*}
    u=\sum_{I\in {\cal I}}\sum_{\gamma\in Q^+(\pi)}\sum_{\beta\in P(\pi)}B_Ia_{I\gamma\beta}
  \end{align*}  
  where $a_{I\gamma\beta}\in G^+_{\gamma}\tau(\beta)$ and the set
  $\{a_{I\gamma\beta}\,|\,\gamma\in Q^+(\pi)$ and $\beta\in P(\pi)\}$ is linearly independent
  for fixed $I\in \cI$. By Proposition  
  \ref{triangIwa}(i) one gets
  \begin{align*}\sum_{\gamma\in Q^+(\pi)}\sum_{\beta\in
    P(\pi)}a_{I,\gamma,\beta}\in G({\le}(\lambda-{\rm wt}(I)))
  \end{align*} 
  for each $I\in {\cal I}$. Hence $a_{I,\gamma,\beta}\in G({\le}(\lambda-{\rm wt}(I)))$ 
  for all $I\in \cI$, $\gamma\in Q^+(\pi)$, and $\beta\in P(\pi)$.
  The multiplicative property (\ref{mult}) and 
  $B_i\in G({\le}\alpha_i)$ ensure that $B_I\in G({\le}{\rm wt}(I))$.
  Hence $B_Ia_{I\gamma\beta}\in G({\le}\lambda)$ for all $I,\gamma,\beta$.
  On the other hand, the definition of $b$ implies that
  \begin{align*}b(u)=\sum_{I\in {\cal I}}\sum_{\tilde\gamma=0}
    \sum_{\tilde\beta=0}B_Ia_{I\gamma\beta}.
  \end{align*}
  Therefore $b(u)\in G({\le}\lambda)$.
  $\blacksquare$
 
\section{Highest weight vectors  in $G^-{\cal M}^+T'_{\Theta}$}\label{HWVinGMT}
  
Let ${\cal V}$ be the vector space defined by 
\begin{align*}
  {\cal V}= G^-{\cal M}^+T'_{\Theta}\cap F_r(\check U)^{U^-}\cap 
  \sum_{\lambda\in 2P^+(\Sigma)}\check U_{-\lambda}.
\end{align*} 
In other words, ${\cal V}$ is the span of all highest weight vectors of spherical 
weight contained in $G^-{\cal M}^+T'_{\Theta}\cap F_r(\check U)$. 
  
The main result of this section finds a basis of ${\cal V}$.
This basis will be used in Section 8 to find a basis of $Z(\check B)$.  
Both bases are indexed by the subset $P_{Z(\Bc)}$ of $ \Ppiplus$ defined by 
(\ref{PBZ-def}) below. Moreover, in subsection 7.3 we obtain information about the weight 
spaces of spherical submodules  of $F_r(\check U)$ generated by highest weight vectors
not contained in $G^-{\cal M}^+T'_{\Theta}$. 
   
\vspace{.5cm}
\noindent
{\bf 7.1 Highest weight vectors for spherical modules.}
The next proposition determines which values $\mu\in P^+(\pi)$ satisfy the following condition:    
$(\adr U)\tau(2\mu)$ contains a highest weight vector of spherical weight that is also an 
element of $G^-{\cal M}^+T'_{\Theta}$. Recall that $w_0$ is the longest element of the 
Weyl group $W$ and that $w_0'$ is the longest element of the 
subgroup $W(\pi_{\Theta})$ of $W$ (see Section 3.2).  
Let $\Ptheta$ be the subset of $\Ppiplus$ defined by 
\begin{align}\label{PBZ-def}
  \Ptheta=\{\mu\in \Ppiplus\,|\,\Theta(\mu)=\mu+w_0\mu-w_0'\mu\}.
\end{align}

\setcounter{proposition}{0}    
\begin{proposition}\label{Yweightvec}
  Let $\mu\in P^+(\pi)$. Then $(\adr U)\tau(2\mu)\cap {\cal V}\neq \{0\}$ if and only if 
  $\mu\in \Ptheta$. Moreover, in this case any nonzero element 
  $Y_\lambda\in(\adr U)\tau(2\mu)\cap {\cal V}$ is a weight vector of spherical
  weight $\lambda=2\mutil$, and $\dim((\adr U)\tau(2\mu)\cap {\cal V})=1$.
\end{proposition}
	       
\noindent
  {\bf Proof:} 
  Fix $\mu\in \Ppiplus$. Assume first that $\mu\in \Ptheta$ and set 
  $\lambda=2\tilde\mu$.  By \cite[Section 2]{a-Letzter-memoirs}, 
  $\{\tilde\gamma|\gamma\in {P^+(\pi)}\}$ is a subset of $P^+(\Sigma)$.
  Hence $\lambda$ is a spherical weight. Note that 
  $2\tilde\mu=\mu-\Theta(\mu)$. Thus, the definition 
  of $\Ptheta$ ensures that $\lambda=w_0'\mu-w_0\mu$. Hence 
  \begin{align}\label{muminusw0mu}
    \mu-w_0\mu-\lambda=\mu-w_0'\mu\in \N\pi_{\Theta}.
  \end{align}
  By Theorem \ref{sphersubprop}, there exists a highest weight 
  vector $Y_{\lambda}\in (\adr U)\tau(2\mu)$ of weight $\lambda$ with respect to the right 
  adjoint action. Moreover, $Y_{\lambda}$ is unique up to nonzero scalar 
  multiplication. Thus the assertion $\dim((\adr U)\tau(2\mu)\cap {\cal V})=1$ follows once 
  we establish below that $Y_{\lambda}\in G^-{\cal M}^+T'_{\Theta}$.
  
  In the graded algebra $\Gr$ the corresponding highest weight vector
  $\oYlam\in (\adr U)\otau$ of weight $\lambda$ satisfies
  \begin{align*}
    \oYlam\in \sum_{\alpha\leq \mu-w_0\mu} K(2\mu)^-_\alpha \ot\cC\otau \ot 
    K(2\mu)^+_{\lambda-\alpha}.
  \end{align*}
  By Lemma \ref{adra} and the fact that $Y_{\lambda}$ has weight 
  $\lambda$ with respect to the right adjoint action, there exists 
  \begin{align*}
    a\in\sum_{\zeta\leq \mu-w_0\mu}
    U^-_{-\zeta}U^+_{\zeta-\lambda}
  \end{align*}
  such that $(\adr a){\overline{\tau(2\mu)}} =\oYlam$.
  Inclusion (\ref{adUzetaxi}) ensures that
  \begin{align}\label{assertion0}
     Y_{\lambda}\in \sum_{\zeta\leq \mu-w_0\mu}\,\,\sum_{\alpha\leq \zeta}\,\,
      \sum_{0\le \eta\le \zeta-\lambda}
     G^-_{-\alpha}U^+_{\alpha-\lambda}\tau(2\mu-\lambda-2\eta).
  \end{align} 
  Note that the inequalities of weights under the summations imply $\alpha-\lambda\in \N\pt$ 
  and $\eta\in \N\pt$. Moreover, by definition $\lambda=\mu-\Theta(\mu)$ and hence
  \begin{align}\label{assertion2}
    \tau(2\mu-\lambda-2\eta)=\tau(\mu+\Theta(\mu)-2\eta)\in T'_{\Theta}.
  \end{align}  
  Combining (\ref{assertion0}), (\ref{assertion2}), and $\alpha-\lambda\in \N\pt$ one 
  obtains $Y_{\lambda}\in  G^-{\cal M}^+T'_{\Theta}$ and hence 
  $Y_\lambda\in {\cal V}\cap (\adr U)\tau(2\mu)$.	  
	  
  Now assume that $(\adr U)\tau(2\mu)$ contains a highest weight 
  spherical vector $Y_{\lambda}$ of weight $\lambda$  with respect to 
  the right adjoint action.  Assume further that
  $Y_{\lambda}\in G^-{\cal M}^+T'_{\Theta}$. 
  In view of the  isomorphism  $g:(\adr U)\tau(2\mu)\rightarrow(\adr U)\overline{\tau(2\mu)}$   given by 
  (\ref{griso}), Lemma \ref{hwcor} implies that there exists a nonzero 
  $Z\in U^-_{-\mu+w_0\mu}G^+\tau(2\mu)$ such that 
  \begin{align*}
    Y_{\lambda}\in Z+\sum_{\alpha\geq 0}\sum_{\gamma<\mu-w_0\mu}
    U^-_{-\gamma}G^+\tau(2\mu-\alpha).
  \end{align*} 
  Since $Y_{\lambda}$ has weight $\lambda$, it follows that $Z$ is a weight 
  vector of weight $\lambda$. Hence 
  \begin{align*}
    Z\in G^-_{-\mu+w_0\mu}U^+_{\mu-w_0\mu-\lambda}\tau(2\mu-\lambda).
  \end{align*}
  The fact that $Y\in G^-{\cal M}^+T'_{\Theta}$ ensures that $Z$ is in 
  $G^-_{-\mu+w_0\mu}{\cal M}^+ T'_{\Theta}$. 
  Thus  $\mu-{w_0\mu}-\lambda \in \N \pi_{\Theta}$ and $\Theta(2\mu-\lambda)=2\mu-\lambda$. 
  Since $\lambda$ is spherical, we have that 
  $\Theta(\lambda)=-\lambda$. Hence  $\lambda=\mu-\Theta(\mu)=2\tilde\mu$.
  Theorem \ref{sphersubprop} further implies that $\lambda=w_0'\mu-w_0\mu$. 
  Thus $\mu-\Theta(\mu)=w_0'\mu-w_0\mu$ and so $\mu\in \Ptheta$.  
  $\blacksquare$
		  
\vspace{.5cm}
\noindent
{\bf 7.2 Highest weight spherical vectors in $G^-{\cal M}^+T'_{\Theta}$.}
Fix $\mu\in \Ptheta$. By Proposition \ref{Yweightvec}, there exists a unique (up to nonzero 
scalar multiplication) highest weight vector of   
weight $2\tilde\mu$ in $G^-{\cal M}^+T'_{\Theta}\cap (\adr U)\tau(2\mu)$ with respect to 
the right adjoint action. Denote this highest weight vector by $Y^{\mu}$. Note that 
$\mu=0$ is in $\Ptheta$ and $Y^{0}$ is just a nonzero scalar. Moreover, 
$Y^{\mu}\in {\cal V}$ for all $\mu\in \Ptheta$.     

\setcounter{theorem}{1}
\begin{theorem}\label{Ybasis} 
  The set of $\{Y^{\mu}|\ \mu\in \Ptheta\}$ is basis for ${\cal V}$.
\end{theorem}
	     
\noindent
{\bf Proof:} 
Note that each $Y^{\mu}\in (\adr U)\tau(2\mu)\setminus\{0\}$.   
Hence the set $\{Y^{\mu}\,|\,\mu\in \Ptheta\}$ is linearly independent.

Let $Y\in {\cal V}$.  Note that ${\cal V}$ can be written as a direct sum of its weight 
spaces. Hence we may reduce to the case where $Y$ is a weight vector of weight $\lambda$. 
Moreover, the definition of ${\cal V}$ ensures that   $\lambda$ is a spherical weight.
Since $Y\in F_r(\check U)$, we can write $Y=\sum_{\gamma\in \Ppiplus}a_{\gamma}$ with 
$a_{\gamma}\in (\adr U)\tau(2 \gamma)$.
It follows from (\ref{adUzetaxi}) that
\begin{align}\label{a-gamma-in}
  a_\gamma\in \sum_{\eta,\xi\ge 0}U^-_{-\lambda-\xi}G^+_\xi\tau(2\gamma-2\eta)\subseteq
  \sum_{\eta\ge 0}G^- U^+\tau(2\gamma-\lambda-2\eta).
\end{align}
By Proposition \ref{Yweightvec} it suffices to show that 
\begin{align}\label{agammaGoal}
  a_\gamma\in G^-\cM^+T_\Theta'
\end{align}  
for all $\gamma$. We prove relation (\ref{agammaGoal}) by induction on 
$m=\max\{h(\gamma)\,|\, a_\gamma\neq 0)\}$. If $m=0$ then $Y=a_0$ and (\ref{agammaGoal})
holds. Assume $m>0$ and let $\oY=\sum_{h(\gamma)=m}\bar{a}_\gamma$
denote the element represented by $Y$ in $\Gr_{2m}$.
The relation $Y\in G^-\cM^+T_\Theta'$ implies
\begin{align*}
  \bar{a}_\gamma\in \sum_{\alpha\in \N\pt} K(2\gamma)^-\ot 
  \overline{\tau(2 \gamma)}\ot K(2\gamma)^+_\alpha
\end{align*}
for all $\gamma$ such that $h(\gamma)=m$. By Lemma \ref{adra} using 
$a_\gamma\in \Uc_{-\lambda}$ one obtains
\begin{align*}
  a_\gamma\in \sum_{\xi\in \N\pt}(\adr U_{-\lambda-\xi}^- U^+_\xi)\tau(2\gamma).
\end{align*}
Hence by (\ref{adUzetaxi}) for all $\gamma$ with $h(\gamma)=m$ one obtains the 
following refinement of (\ref{a-gamma-in})
\begin{align}\label{a-gamma-in2}
  a_\gamma\in \sum_{\eta,\xi\in \N\pt}U^-_{-\lambda-\xi}G^+_\xi\tau(2\gamma-2\eta)\subseteq
    \sum_{\eta\in \N\pt}G^- \cM^+\tau(2\gamma-\lambda-2\eta). 
\end{align}
For all $\gamma$ with $h(\gamma)=m$, as $\bar{a}_\gamma\neq 0$, there exists a nonzero 
\begin{align*}
  X_\gamma\in \sum_{\xi\in \N\pt}U^-_{-\lambda-\xi}G^+_\xi\tau(2\gamma)\subseteq
  G^- \cM^+\tau(2\gamma-\lambda) 
\end{align*}
such that 
\begin{align*}
  a_\gamma\in X_\gamma+\sum_{\eta<\gamma}\sum_{\xi\in \N\pt} 
  U^-_{-\lambda-\xi}G^+_\xi\tau(2\eta)\subseteq X_\gamma +\sum_{\eta<\gamma} G^-\cM^+
  \tau(2\eta-\lambda).
\end{align*}
The maximality of $\gamma$ and the relation 
$Y\in G^-\cM^+T_\Theta'$ now imply that $\tau(2\gamma-\lambda)\in T_\Theta'$ for all $\gamma$
such that $h(\gamma)=m$. Together with (\ref{a-gamma-in2}) this yields 
$a_\gamma\in G^-\cM^+T_\Theta'$ for all $\gamma$ such that $h(\gamma)=m$. Hence
\begin{align*}
  \sum_{h(\gamma)<m}a_\gamma=Y-\sum_{h(\gamma)=m}a_\gamma\in G^-\cM^+T_\Theta'
\end{align*}
and by the inductive hypothesis $a_\gamma\in G^-\cM^+T_\Theta'$ for all $\gamma$.
$\blacksquare$			 

\vspace{.5cm}
\noindent
{\bf 7.3 $B$-invariant subspaces of $G({\le}2\mu)$.}
Fix $\mu\in \Ptheta$ and recall the definition of  $Y^{\mu}$  given in 
Section 7.2. Recall (\ref{vBweights}) and let $v_{\mu}^B$ 
be the  uniquely determined nonzero $(\adr B)$-invariant element of $(\adr U)Y^{\mu}$ such that
\begin{align}\label{vmuBdefn}
  v_{\mu}^B\in Y^{\mu}+\sum_{\{\gamma\in P(\Sigma)|\ \gamma<_r\lambda\}}\check U_{-\gamma}.
\end{align}
Recall further the definition of $G({\le}2\mu)$ given in (\ref{Glessthanmu}).  
Set $V_0^{\mu}$ equal to  the sum of all the spherical submodules of $ G({\le}2\mu)$ 
of highest weight unequal to $2\tilde\mu$. 
      
Consider $\lambda\in P^+(\pi)$.  As discussed at the beginning of Section 3.1, 
the highest weight of $(\adr U^-)\overline{\tau(2\lambda)}$ with respect to the right adjoint 
action is $\lambda-w_0\lambda$.   It follows from the 
description of $(\adr U)\overline{\tau(2\lambda)}$ in (\ref{grinto}) 
that the highest weight of $(\adr U)\overline{\tau(2\lambda)}$ is 
$\lambda-w_0\lambda$. Using the isomorphism (\ref{griso}) we see that 
\begin{align}\label{adrUproperty}
  (\adr U)\tau(2\lambda)\subseteq \sum_{\gamma\leq \lambda-
  w_0\lambda}\check U_{-\gamma}.
\end{align}

\setcounter{lemma}{2}      
\begin{lemma}\label{G2muBdecomp} 
  Let $\mu\in \Ptheta$. Then $((\adr U)\tau(2\eta))^B\subseteq (V_0^{\mu})^B$ 
  for all $\eta\in P^+(\pi)$ satisfying $\eta<\mu$ and $\tilde\eta<\tilde\mu$.  
  Moreover, $(V_0^{\mu})^B= G({\le}2\mu)^B\cap \sum_{ h^{\Sigma}(\gamma)< h^{\Sigma}(2\mu)}
  \check U_{-\gamma}$.
\end{lemma} 

\noindent
{\bf Proof:} The assumption $\eta<\mu$ forces $-w_0\eta<-w_0\mu$ and hence 
$\eta-w_0\eta<\mu-w_0\mu$. Since $\etatil<_r\mutil$ we further get 
$\tilde\eta-\widetilde{w_0\eta}<_r \tilde\mu-\widetilde{w_0\mu}$.
Hence by 
(\ref{adrUproperty}) and (\ref{vBweights}), we see that 
\begin{align}\label{adrutau2etas}
  ((\adr U)\tau(2\eta))^B\subseteq \sum_{\gamma<_r \tilde{\mu}-\widetilde{w_0\mu}}\check 
  U_{-\gamma}.
\end{align} 
Note that $\widetilde{w_0'\mu}=\tilde\mu$ and hence $-\widetilde{w_0\mu}=\tilde\mu$ for all 
$\mu\in \Ptheta$. Suppose that $\beta$ is a spherical weight such that 
$V(\beta)^*$ is isomorphic to a submodule of $(\adr U)\tau(2\eta)$.  
By (\ref{adrutau2etas}), we see that 
$\beta<_r\tilde{\mu}-\widetilde{w_0\mu}=2\tilde\mu$.   Hence 
$\beta\neq 2\tilde\mu$. Therefore, $((\adr U)\tau(2\eta))^B\subseteq (V_0^{\mu})^B$ which 
proves the first assertion. To prove the second assertion suppose first that $V$ 
is a finite-dimensional simple spherical submodule of $G({\le}2\mu)$ contained in
$\sum_{h^{\Sigma}(\gamma)<h^{\Sigma}(2\mu)}\check U_{-\gamma}$.  It
follows that the highest weight of $V$ is unequal to $2\tilde\mu$. Hence
\begin{align*}
  G({\le}2\mu)^B\cap(\sum_{ h^{\Sigma}(\gamma)< h^{\Sigma}(2\mu)}
  \Uc_{-\gamma})\, \subseteq (V_0^{\mu})^B.
\end{align*}

It now remains to show that
\begin{align*}
  (V_0^{\mu})^B\cap (\adr U)\tau(2\eta)\subseteq\sum_{h^{\Sigma}(\gamma)<h^{\Sigma}(2\mu)}
  \Uc_{-\gamma}
\end{align*}
for all $\eta\le\mu$.
Suppose that  $V(\beta)^*$ is isomorphic to a submodule of $(V_0^{\mu})^B$ with $\beta$ a 
spherical weight. It follows that $\beta=\tilde\beta$. The highest weight space of 
$(\adr U){\tau(2\eta)}$ is $\eta-w_0\eta\le\mu-w_0\mu$.
Hence $\beta\le_r \etatil - \widetilde{w_0\eta}\le_r \mutil - \widetilde{w_0\mu}=2\mutil$.
The definition of $V_{\mu}^0$ forces $\beta\neq 2\tilde\mu$ and thus $\beta <_r 2\mutil$.
$\blacksquare$

\vspace{.5cm}
   
\section{A basis for the center $Z(\check B)$}\label{ZBbasis}
Fix $\mu\in \Ptheta$ and recall the definitions of $Y^{\mu}$, 
$v_{\mu}^B$, and $V^{\mu}_0$ from Sections 7.2 and 7.3. In this section we show how 
to use the spherical vector 
$v_{\mu}^B$  to construct an element $d_{\mu}$  in $Z(\check B)\cap 
G({\le}2\mu)$ with a nonzero component in $(\adr U)\tau(2\mu)$. We then prove that the set 
$\{d_{\mu}|\ \mu\in \Ptheta\}$ is a basis for $Z(\check B)$.
  
\vspace{.5cm}
\noindent
{\bf 8.1 Formulas for the right adjoint action.}
  Let $\alpha_k\in \pi\setminus\pi_{\Theta}$.  By \cite[(7.15)]{MSRI-Letzter} we have 
the following formula for the coproduct of $B_k$ 
\begin{align}\label{kowBkformula}
  \kow (B_k)\in t_k\ot B_k+U\ot \cM\Tt. 
\end{align}
Recall moreover from (\ref{Bi-def}) that $\vep(B_k)=s_k$ for all $\alpha_k\in\pi\setminus\pt$.
\setcounter{lemma}{0}
\begin{lemma} \label{adBG} 
  If $G$ is an $(\adr \cM\Tt)$-invariant element of 
  $\check U$, then
  \begin{align*}
    (\adr B_k)G-s_kG=t_k^{-1}(GB_k-B_kG) 
  \end{align*} 
  for all $\alpha_k\in \pi\setminus\pi_{\Theta}$.
\end{lemma}
 
\noindent{\bf Proof:}  
Fix $\alpha_k\in \pi\setminus\pi_{\Theta}$. By (\ref{kowBkformula}) we have 
${\it\Delta} (B_k)=t_k\otimes B_k +\sum_ia_i\otimes b_i$
where $a_i\in U$ and $b_i\in {\cal M}T_{\Theta}$. Hence
$s_k=\epsilon(B_k)=\sigma(t_k)B_k+\sum_i\sigma(a_i)b_i$ and therefore 
$\sum_i\sigma(a_i)b_i=-t_k^{-1}B_k+s_k$. Let $G$ be an $(\adr {\cal M}T_{\Theta})$-invariant 
element of $\check U$. Since ${\cal M}T_{\Theta}$ is a Hopf subalgebra of $\check U$, it 
follows that $G$ commutes with elements of ${\cal M}T_{\Theta}$.  Hence
\begin{align*}
  (\adr B_k)G=&t_k^{-1}GB_k+\sum_i\sigma(a_i)Gb_i 
    =t_k^{-1}GB_k+\sum_i\sigma(a_i)b_iG\\
    =&t_k^{-1}GB_k-t_k^{-1}B_kG + s_kG.\, \blacksquare
\end{align*} 
    
\vspace{.5cm}
    
\noindent{\bf 8.2 $B$-invariant elements of $G({\le}2\mu)$.}
The next technical lemma will be used to complete the inductive 
step in  the proof of Theorem \ref{mainThm}.
\setcounter{lemma}{1}     
\begin{lemma}\label{mainlemma} 
  Let $\mu\in \Ptheta$ and $v\in v_{\mu}^B+ (V_0^{\mu})^B$.
  Assume that there exists $m\in {{1}\over{N}}\Z$ such that
  \begin{align}\label{v-bv-smallweights}
    v-b(v)\in \sum_{h^{\Sigma}(\eta)<m}\check U_{-\eta}.
  \end{align}
  Then there exists $w\in v_{\mu}^B+(V_0^\mu)^B$ such that 
  $w-b(w)\in  {\cal F}^{\theta}_{m-{{1}\over{N}}}(\check U)$.
\end{lemma}

\noindent {\bf Proof:} 
Without loss of generality we may assume that
\begin{align}\label{v-bv-largefilter}
  v-b(v)\not\in {\cal F}^{\theta}_{m-{{1}\over{N}}}(\check U)
\end{align}
because otherwise $w=v$ fulfills the claim of the lemma.
   
Since $v\in \check U^B$, it follows that $v$ is $(\adr {\cal M}T_{\Theta})$-invariant.
Recall that $b$ is a projection of ${\cal M}T_{\Theta}$-modules.
Hence $b(v)$ is also $(\adr {\cal M}T_{\Theta})$-invariant.   
Fix $\alpha_k\in\pi\setminus\pi_{\Theta}$. By Lemma \ref{adBG} we have
\begin{align*} 
  (\adr B_k)(v-b(v))-s_k(v-b(v))=t_k^{-1}[(v-b(v))B_k-B_k(v-b(v))].
\end{align*}
Hence, (\ref{v-bv-smallweights}) and the definition (\ref{Bi-def}) of $B_i$ ensure that  
\begin{align}\label{adrBkvB}  
  (\adr B_k)(v-b(v))-s_k(v-b(v))\in \sum_{h^{\Sigma}(\eta)< m+1} \check U_{-\eta}.
\end{align} 
On the other hand, the fact that $v\in \check U^B$ and $b(v)\in \check B$ 
implies that   
\begin{align*} 
  (\adr B_k)(v-b(v))-s_k(v-b(v)) &=t_k^{-1}[b(v)B_k-B_kb(v)]\\
    &\in t_k^{-1}\sum_{I\in {\cal I}}B_I{\cal M}^+T'_{\Theta}.
\end{align*}
Let $r$ be the smallest integer such that 
\begin{align}\label{tk} 
  (\adr B_k)(v-b(v))-s_k(v-b(v))= t_k^{-1}\sum_{h^{\Sigma}({\rm wt}(I))\leq r}B_Im_I
\end{align}
where each $m_I\in {\cal M}^+T'_{\Theta}$.  Note that the assumption on $r$ forces at least one 
element of the set $\{m_I|\ h^{\Sigma}({\rm wt}(I))=r\}$ to be nonzero.
By Lemma \ref{BIdecomp}, it follows that
\begin{align}\label{tk2} 
  (\adr B_k)(v& -b(v))-s_k(v-b(v))\nonumber\\
  &\in
  \sum_{h^{\Sigma}({\rm wt}(I))= r}t_k^{-1}\Big[(yt)_Im_I{+}
  \sum_{\eta<{\rm wt}(I)}U^-_{-\eta}G^+\check U^0\Big].
\end{align}
Examining the right hand side of (\ref{tk2}), we see that the lowest weight
term of $(\adr B_k)(v-b(v))-s_k(v-b(v))$ written as a sum of nonzero weight
vectors is contained in $\sum_{h^{\Sigma}(\eta)=r}\check U_{-\eta}$.  	   
Hence (\ref {adrBkvB}) implies that $r< m+1$. Note moreover that
$B_I\in \cF^\theta_r(\Uc)$ if $\hs(\wght(I))\le r$.
It now follows directly from (\ref{tk}) that
\begin{align}\label{last}
  (\adr B_k)(v-b(v))-s_k(v-b(v))\in {\cal F}^{\theta}_{m-\frac{1}{N}}(\check U).
\end{align}
Let $s$ denote the degree of $v-b(v)$ with respect to the filtration ${\cal F}^{\theta}$.
By (\ref{v-bv-largefilter}), we see that $s\geq m$. By definition of $B_k$ an 
$(\cM\Tt)$-invariant element $w$ of a $B$-module is $B$-invariant if and only if
$B_k w-s_kw=0$ for all $\alpha_k\in \pi\setminus \pt$.
Hence (\ref{last}) implies that the 
image of 
$v-b(v)$ in ${\cal F}^{\theta}_s(\check U)/{\cal F}^{\theta}_{s-{{1}\over{N}}}(\check U)$ 
is $B$-invariant.
 
By assumption, $v\in G({\le}2\mu)$.  Hence, by Proposition \ref{constructingLem2}(i)
one has $b(v)\in G({\le}2\mu)$. It follows that $h^{\Sigma}(2\mu)\geq s\geq m$.
Theorem \ref{GrB} and Proposition \ref{FGrtheta} ensure  that there exists 
$v'\in G({\le}2\mu)^B$ such that
\begin{align*}
  v-b(v)\in v'+{\cal F}^{\theta}_{m-{{1}\over{N}}}(\check U).
\end{align*}
The relation $h^{\Sigma}(2\mu)\geq m$, assumption (\ref{v-bv-smallweights}), and Lemma 
\ref{G2muBdecomp} imply that $v'\in (V_0^{\mu})^B$.
 
Set $w=v-v'$.  Note that $w\in v+(V_0^{\mu})^B=v_{\mu}^B+(V_0^{\mu})^B$.
By the definition of the projection map $b$, we have that $v-b(v)\in  B\Ttc {\cal G}$.
Hence $v'\in B\Ttc {\cal G}+{\cal F}^{\theta}_{m-{{1}\over{N}}}(\check U)$.
By Proposition \ref{constructingLem2}(ii), we obtain
$b(v')\in {\cal F}^{\theta}_{m-{{1}\over{N}}}(\check U)$.
Therefore $w-b(w)=[v-b(v)-v']+b(v')\in {\cal F}^{\theta}_{m-{{1}\over{N}}}(\check U).$
$\blacksquare$     
 
\vspace{.5cm}
\noindent{\bf 8.3 Finding central elements.}
Suppose that $\mu\in P^+(\pi)$. The next theorem shows how 
to construct an  element of $Z(\check B)$ that is also contained in $G({\le} 2\mu)$. 
\setcounter{theorem}{2}
\begin{theorem}\label{mainThm}
  Given $\mu\in \Ptheta$, there exists an 
  element $d_{\mu}\in Z(\check B)$ which lies in $v_{\mu}^B+(V_0^{\mu})^B$.
  Moreover, 
  \begin{align}\label{weightdmu}
    d_{\mu}\in Y^{\mu}+\sum_{h^{\Sigma}(\gamma)<h^{\Sigma}(2\mu)}\check U_{-\gamma}.
  \end{align}
\end{theorem}
{\bf Proof:} 
  Set $V=(V_0^{\mu})^B$. The definition (\ref {vmuBdefn}) of $v_{\mu}^B$ and Lemma
  \ref{G2muBdecomp} imply that 
  $v_{\mu}^B+V\subseteq Y^{\mu}+\sum_{h^{\Sigma}(\gamma)<h^{\Sigma}(2\mu)}\Uc_{-\gamma}$.  
  Hence (\ref{weightdmu}) follows from the first assertion of the theorem.
  
  Note that the first assertion of the theorem holds provided we prove
  $(v_{\mu}^B+V)\cap \check B\neq \{0\}$.  It suffices to show that there exists 
  $v\in v_{\mu}^B+V$ such that $v=b(v)$. Assume on the contrary that $v-b(v)\neq 0$ for all 
  $v\in v_{\mu}^B+V$. It follows that for each $v\in v_{\mu}^B+V$, there exists a smallest 
  rational number $n(v)$ such that 
  \begin{align*}
    v-b(v)\in {\cal F}_{n(v)}^{\theta}(\Uc).
  \end{align*}
  Recall that $V$ is finite dimensional. Hence the set $\{n(v)\,|\, v\in v_{\mu}^B+V\}$
  has a minimal element $n_{min}$. In what follows we will construct an element 
  $w\in v_{\mu}^B+V$ such that $w-b(w)\in {\cal F}_{n_{min}-1/N}^{\theta}(\Uc)$ leading
  our assumption $(v_{\mu}^B+V)\cap \check B\neq \{0\}$ to a contradiction.

  \medskip
  \noindent
  {\bf Claim 1:} $n_{min}<h^{\Sigma}(2\mu)$.
     
  \medskip
  \noindent
  \textit {Proof of Claim 1:} Recall that $Y^{\mu}\in G^-{\cal M}^+T'_{\Theta}$ is a 
  vector of weight $2\tilde\mu$ with respect to the right adjoint action. Hence, we can write
  \begin{align}\label{Ylambdaone}
    Y^{\mu}\in\sum_{\tilde\beta=2\tilde\mu}G^-_{\beta}{\cal M}^+T'_{\Theta}.
  \end{align}
  Using Lemma \ref{BIdecomp} and induction one obtains
  \begin{align*}
    Y^{\mu}&\in\sum_{h^{\Sigma}({\rm wt}(I))=\hs(2\mu)}
    B_I{\cal M} T'_{\Theta}+\sum_{h^{\Sigma}({\rm wt}(J)) <\hs(2\mu)}B_JG^+\Uc^0\\
    &\subseteq \sum_{h^{\Sigma}({\rm wt}(I))=\hs(2\mu)}
    B_I{\cal M} T'_{\Theta}+\sum_{h^{\Sigma}(\eta) <\hs(2\mu)}U^-_{-\eta}G^+\Uc^0.
  \end{align*}
  It follows from the definition  (\ref{vmuBdefn}) of $v_{\mu}^B$ that
  \begin{align}\label{v0insum}
    v_{\mu}^B\in \sum_{h^{\Sigma}({\rm wt}(I))=\hs(2\mu)}B_I{\cal M} 
    T'_{\Theta}+\sum_{h^{\Sigma}(\eta)< \hs(2\mu)}\check U_{-\eta}.
  \end{align}
  By Proposition \ref{constructingLem2}(iii), applying the projection $b$ to both sides of
  (\ref{v0insum}) yields
  \begin{align}\label{vBminusbvB}
    v_{\mu}^B-b(v_{\mu}^B)\in \sum_{h^{\Sigma}(\eta)<\hs(2\mu)} \Uc_{-\eta}.
  \end{align}
  Hence we can apply Lemma \ref{mainlemma} to $v_{\mu}^B$ and there
  exists $w\in v_{\mu}^B+V$ with $n(w)<\hs(2\mu)$.  Therefore, $n_{\min}<\hs(2\mu)$. 
  $\square$

  \medskip
  We now continue with the proof of Theorem \ref{mainThm}. Consider $v\in v_{\mu}^B+V$ such that
  $n(v)=n_{min}$. By definition of $V$ we can write $v=v_{\mu}^B +v'+v''$ where
  \begin{align*}
    v'\in \sum_{\makebox[0cm]{$\eta\leq \mu \atop h^{\Sigma}(2\eta)> n_{min}$}}
	   ((\adr U)\tau(2\eta))^B,\qquad
    v''\in \sum_{\makebox[0cm]{$\eta\leq \mu \atop h^{\Sigma}(2\eta)\leq n_{min}$}}
	 ((\adr U)\tau(2\eta))^B.
  \end{align*}
  It follows from Proposition \ref{constructingLem2} (ii) that 
  $v''-b(v'')\in {\cF}^{\theta}_{n_{min}}(\check U)$.  Hence subtracting $v''$ we may assume 
  that $v=v_{\mu}^B +v'$ and $n(v)=n_{min}$.
	    
  By Lemma \ref{BtrigLemone}, we see that $b(v)\in G^-U^+ T_{\geq}$.
  It follows from Corollary \ref{v-bcor} that $v\in G^-U^+ T_{\geq}$.
  The definition of the degree function defining ${\cal F}^{\theta}$ implies that
 \begin{align} 
   v-b(v)&\in \sum_{\eta\geq 0}\sum_{\beta\geq 0}\sum_{\tilde\gamma\geq_r 0}
   \sum_{h^{\Sigma}(\eta+\beta+\gamma)\le n_{min}}(G^-_{-\eta}U_{\beta}^+)\tau(\gamma)
   \nonumber\\
   &\subseteq \sum_{h^{\Sigma}(\eta)=n_{min}}G^-_{-\eta}{\cal M}^+T'_{\Theta}+
   \sum_{h^{\Sigma}(\eta)<n_{min}}(G^-_{-\eta}U^+)T_{\geq}.\label{vb(v)subset}
 \end{align}
 Applying  Lemma \ref{BIdecomp} and induction to (\ref{vb(v)subset}) yields
 \begin{align}\label{anotherlabel}v-b(v)&\in
  \sum_{h^{\Sigma}({\rm wt}(I))=n_{min}}B_I{\cal M}^+T'_{\Theta}+
  \sum_{h^{\Sigma}({\rm wt}(I))<n_{min}}B_IU^+\check U^0.
 \end{align}
 By the definition of the projection map $b$, we have that
 $v-b(v)\in B\check T_{\Theta}{\cal G}$.   Hence the definition of
 ${\cal G}$  in Section 6.1 combined with (\ref{anotherlabel}) imply  that
 \begin{align*}
   v-b(v)&\in\sum_{h^{\Sigma}({\rm wt}(I))<n_{min}}B_IU^+\check U^0.
 \end{align*}
 Hence by Lemma \ref{BIdecomp}, we see that
 \begin{align}\label{noidea}
   v-b(v)\in \sum_{h^{\Sigma}(\eta)<\nm}G^-_{-\eta}U^+\check
   U^0\subseteq \sum_{h^{\Sigma}(\eta)<\nm}\check U_{-\eta}.
 \end{align}
 Assertion (\ref{noidea}) allows us to apply Lemma \ref{mainlemma} to $v$.
 Thus there exists $w\in v_{\mu}^B+V$ such that $n(w)\leq n_{min}-{{1}\over{N}}$.
 This is the desired contradiction to the minimality of $n_{min}$.
 $\blacksquare$

\vspace{.5cm}
\noindent
{\bf 8.4 Central elements in $(\adr U)\tau(2\mu)$.}
We show here that if $\mu\in \Ptheta$ satisfies a suitable minimality
condition then Theorem \ref{mainThm} yields elements in 
$Z(\Bc)\cap(\adr U)\tau(2\mu)$.
These elements are used in \cite{a-Ktobe} to obtain solutions of the reflection equation
and hence further establish the link between different approaches to the construction of
quantum symmetric pairs.
\setcounter{corollary}{3}
\begin{corollary}
  Fix $\mu\in \Ptheta$. Assume that for all $\nu\in \Ppi$ the 
  condition $\nu<\mu$ implies $\nu\notin \Ppiplus\setminus\{0\}$.
  Then there exists a nonzero element $d_\mu\in Z(\Bc)\cap (\adr U)\tau(2\mu)$. 
\end{corollary}
\noindent{\bf Proof:} The definition of $(V^\mu_0)^B$ in Section 7.3 and 
Theorem \ref{locfinBthm} imply 
\begin{align*}
  (V^\mu_0)^B \subseteq \sum_{\nu\in \Ppiplus,\,\nu\le\mu}(\adr U)\tau(2\nu).
\end{align*}  
The assumption about $\mu$ now implies that
$(V^\mu_0)^B$ is a subspace of $(\adr U)\tau(2\mu)+(\adr U)\tau(0)$.
$\blacksquare$

\vspace{.5cm}
\noindent
{\bf 8.5 A basis for the center of $\check B$.}
 By Theorem \ref{mainThm}, we can choose an element $d_\mu\in Z(\check B)$ which lies
 in $v_{\mu}^B+(V_0^{\mu})^B$ for each $\mu\in \Ptheta$.  We will assume these elements 
 to be fixed for the rest of this paper. Note that when $\mu=0$, then 
 $d_0=v_0^B=Y_0$ is a nonzero scalar.
  \setcounter{theorem}{4}
 \begin{theorem}\label{basis-prop}
   The set $\{d_\mu\,|\,\mu \in \Ptheta\}$ is a basis of the vector space $Z(\Bc)$.
 \end{theorem}
 \noindent{\bf Proof:} 
   Note that $d_{\mu}$ lies in $(\adr U)\tau(2\mu)\setminus\{0\}$ up to 
   terms of lower filter degree with respect to ${\cal F}$.   
   Hence the set $\{d_{\mu}|\mu\in \Ptheta\}$ is   linearly independent.
    
   Let $C$ in $Z(\check B)$.  Write $C=\sum_{\beta}c_\beta$ with 
   $c_\beta\in \check U_{-\beta}$. Set 
   \begin{align*}
     m={\rm max}\{h^{\Sigma}(\beta)|\ c_{\beta}\neq 0\}.
   \end{align*} Let $\{\beta_i|\  1\leq i\leq k\}$ denote the set of  weights such 
   that $h^{\Sigma}(\beta_i)=m$ and $c_{\beta_i}\neq 0$.  
   Since $C$ is a sum of spherical vectors inside $F_r(\check U)$, 
   it follows from (\ref{vBweights}) that each $c_{\beta_i}$ is a  highest 
   weight  vector with respect to the right adjoint action contained in a 
   finite-dimensional spherical submodule of 
   $F_r(\check U)$.   
   By (\ref{Bdirectsum}) and the fact that $C\in\check B$, we can write 
   \begin{align}\label{C-triang}
     C\in\sum_{I\in {\cal I}} B_I {\cal M}^+T'_{\Theta}.   
   \end{align} 
   It follows from Lemma \ref{BIdecomp}, that $c_{\beta_i}  \in G^-{\cal M}^+T'_{\Theta}$ 
   for $i=1,\dots, k$. By Theorem \ref{Ybasis}, there exist scalars $a_{i\gamma}$ such that 
   \begin{align*}c_{\beta_i}=
     \sum_{\gamma\in \Ptheta,\ h^{\Sigma}(\gamma)=m}a_{i\gamma}Y_{\gamma}.
   \end{align*} Note that the fact that each $\beta_i$ is a 
   spherical weight forces $m=h^{\Sigma}(\beta_i)\geq 0$.
   If $m=0$, then $i=1$, $\beta_1=0$, and  $c_{0}$ is just a 
   scalar multiple of the nonzero scalar $Y_0$.  Moreover,   
   the simple $(\adr U)$-module generated by $c_0$ is just ${\cal C}={\cal C}d_0$, so 
   $C\in {\cal C}d_0$. Assume $m>0$.   
   By Theorem \ref{mainThm}, 
   \begin{align*}
     C-\sum_{i=1}^k\sum_{h^{\Sigma}(\gamma)=m}a_{i\gamma}d_{\gamma}
          \quad	\in \quad Z(\check B)\cap \sum_{h^{\Sigma}(\beta)<m}\check U_{-\beta}.
   \end{align*}	
   The theorem follows using induction on $m$. $\blacksquare$
          
   \vspace{.5cm}
   
   \noindent
   
\section{Realizing $Z(\Bc)$ as a polynomial ring}\label{ZBpoly}
We now calculate the set $\Ptheta$ explicitly. It turns out that $\Ptheta$ is closely 
related to the rank of the invariant Lie algebra $\gfrak^\theta$. 
This is used to show that $Z(\Bc)$ is a polynomial ring in $\rk(\gfrak^\theta)$ variables.
To this end it is essential that 
central elements in $G({\le}2\mu)$ for $\mu\in \Ptheta$ are unique up to elements of lower 
filter degree.
  
\vspace{.5cm}
\noindent {\bf 9.1 Determining $\Ptheta$.}
We restrict to the case of irreducible symmetric pairs. 
By \cite[2.5 and 5.1]{a-Araki62}, the pair $(\gfrak,\gfrak^\theta)$ is irreducible if and only
if $\gfrak$ is simple, or $\gfrak=\gfrak'\oplus\gfrak'$ for some simple $\gfrak'$ and 
$\theta$ interchanges the two isomorphic components.

In the second case the corresponding quantum symmetric pair $\Uc$, $\Bc$ is given explicitly
in \cite[Section 7]{a-Letzter03} and it is straightforward to check that 
$\Bc\cong \Uc_q(\gfrak')$ as $\cC$-algebras.
Hence the center $Z(\Bc)\cong Z(\Uc_q(\gfrak'))$ is well known by \cite{a-JoLet2}. 
Here we obtain the same result,
it should however be borne in mind that our calculations heavily depended on the previous
knowledge of the structure of $F_r(\Uc)$. As $\pt$ is empty one gets 
$\Ptheta=\{\mu\in \Ppiplus\,|\, \Theta(\mu)=-w_0\mu\}.$ Let $P^+(\pi')$ denote the dominant
integral weights of $\gfrak'$ with respect to a set of simple roots $\pi'$ for $\gfrak'$.
Without loss of generality we can
assume that $\Ppiplus=P^+(\pi')\times P^+(\pi')$. Then
\begin{align*}
  \Ptheta=\{(\mu,-w_0^1\mu)\,|\, \mu\in P^+(\pi')\}
\end{align*}
where $w_0^1$ denotes the longest element in the Weyl group of $\gfrak'$. 

Recall Araki's list \cite[5.11]{a-Araki62} of all symmetric pairs
$(\gfrak,\gfrak^\theta)$ for simple $\gfrak$. Following \cite[Section 7]{a-Letzter03} the 
parameter $p$ occurring in Araki's list will here be denoted by $r$, and as before 
$n=\rk(\gfrak)$.
\setcounter{proposition}{0}
\begin{proposition} \label{PBZ}
  Assume that $(\gfrak,\gfrak^\theta)$ is an irreducible symmetric pair.
  The subset $\Ptheta$ of $ \Ppiplus$ is determined by the following list.
\begin{enumerate}     
  \item[1)] If $\gfrak=\gfrak'\oplus\gfrak'$ and $\Ppiplus=P^+(\pi')\times P^+(\pi')$
            then $\Ptheta=\{(\mu,-w_0^1\mu)\,|\, \mu\in P^+(\pi')\}$.
  
  \item[2)] If $\gfrak$ is one of $B_n$, $C_n$, $E_7$, $E_8$, $F_4$, or $G_2$
        then $\Ptheta=\Ppiplus$.
        
  \item[3)] If $(\gfrak,\gfrak^\theta)$ is of type $AI$, $AII$, ($DI$, case 3), 
        ($DIII$, case 1), $EI$ , 
        or $EIV$ then $\Ptheta=\{\mu\in\Ppiplus\,|\,\mu=-w_0\mu\}$.  
        
  \item[4)] If $(\gfrak,\gfrak^\theta)$ is of type $AIII$, $AIV$, $EII$, $EIII$,
        or ($DIII$, case 2) then $\Ptheta=\Ppiplus$.
          
  \item[5)] If $(\gfrak,\gfrak^\theta)$ is of type ($DI$, case 1) or $DII$ then
        \begin{align*}
          \Ptheta=
            \begin{cases}
              \Ppiplus&\mbox{if $r$ is even,}\\
              \sum_{i=1}^{n-2}\N \omega_i+\N(\omega_{n-1}+\omega_n)&
              \mbox{if $r$ is odd.}
            \end{cases}
        \end{align*}
          
  \item[6)] If $(\gfrak,\gfrak^\theta)$ is of type ($DI$, case 2) then
        \begin{align*}
          \Ptheta=
            \begin{cases}
              \Ppiplus&\mbox{if $n$ is odd,}\\
              \sum_{i=1}^{n-2}\N\omega_i+\N(\omega_{n-1}+\omega_n)&
                \mbox{if $n$ is even.}
            \end{cases}
        \end{align*} 
\end{enumerate}  
\end{proposition}
{\bf Proof:}  
The diagonal case 1) has already been handled above. Consider now the case where $\gfrak$
is simple. Recall \cite[Section 7]{MSRI-Letzter} that $\Theta(\mu)=-w_0'd(\mu)$ where 
$d$ is a diagram automorphism. Hence, if $d=\id$ then $\Ptheta$ consists 
of all selfdual dominant integral
weights. This proves statements 2) and 3) of the proposition.

Note that for all cases listed in 4) as well as ($DI$, case $1$, $r$ even) and 
($DI$, case $2$, $n$ odd) one has $d=-w_0$. Note moreover that 
$\alpha_i-d(\alpha_i)$ is invariant under $\Theta$ for any $\alpha_i\in \pi$. 
Hence one obtains 
\begin{align*}
  \Theta(\alpha_i)-\alpha_i=\Theta(d(\alpha_i))-d(\alpha_i)=-w_0'\alpha_i+w_0\alpha_i
\end{align*}
for any $\alpha_i\in \pi$. This proves 4), 5) for $r$ even, and 6) for $n$ odd.

The claim of 5) follows from $d(\omega_i)=-w_0'\omega_i=\omega_i$ for 
$i=1,\dots,n-2$ and
\begin{align*}
  d(\omega_{n-1})&=\begin{cases}  \omega_{n-1}&\mbox{if $n-r$ even,}\\
                 \omega_n&\mbox{if $n-r$ odd,}\end{cases}\\
  d(\omega_{n})&=\begin{cases}  \omega_{n}&\mbox{if $n-r$ even,}\\
                 \omega_{n-1}&\mbox{if $n-r$ odd,}\end{cases}
\end{align*}
\begin{align*}
  -w_0'(\omega_{n{-}1})&=\begin{cases}\omega_{n{-}1}{-}\omega_r&\mbox{if $n-r$ even,}\\
                               \omega_{n}{-}\omega_r&\mbox{if $n-r$ odd,}\end{cases}\\
  -w_0'(\omega_{n})&=\begin{cases}\omega_{n}{-}\omega_r&\mbox{if $n-r$ even,}\\
                               \omega_{n{-}1}{-}\omega_r&\mbox{if $n-r$ odd.}\end{cases}             
\end{align*}

Note that for the symmetric pair of type ($DI$, case 2) one has $w_0'=\id$ and hence 
$\mu\in \Ptheta$ if and only if $d(\mu)=-w_0\mu$. This condition holds precisely for all 
dominant integral weights given in 6). $\blacksquare$   
   
\vspace{.5cm}
\noindent{\bf 9.2 The rank of $\gfrak^\theta$.}   
\noindent Note that for every $(\gfrak,\gfrak^\theta)$ there exist 
(up to ordering) uniquely determined elements $\nu_1,\dots,\nu_m\in \Ppiplus$ such that
\begin{align}\label{PBZ-free}
  \Ptheta=\bigoplus_{i=1}^m \N \nu_i.
\end{align}
We define $\rk(\Ptheta):=m$. Moreover, 
let $\rk(\gfrak^\theta)$ denote the rank of the reductive Lie algebra $\gfrak^\theta$.
\setcounter{proposition}{1}
\begin{proposition}\label{rank=rank}
  The relation {\upshape $\rk(\Ptheta)=\rk(\gfrak^\theta)$} holds.
\end{proposition}
{\bf Proof:}
  We may restrict to the case where $\gfrak$ is simple.
  By \cite[Chapter X, Thm.~5.15 (ii)]{b-Helga78} for every $(\gfrak,\gfrak^\theta)$ there 
  exists an automorphism $\nu$ of the Dynkin diagram of $\gfrak$ such that
  $\rk(\gfrak^\theta)$ coincides with the dimension of the fixed point set
  $\hfrak^\nu$ in the corresponding Cartan subalgebra $\hfrak$ of $ \gfrak$.
  The order $k$ of $\nu$ can be read off \cite[p.~514, Tables II and III]{b-Helga78}.
  In particular if $\gfrak$ has no nontrivial diagram automorphism one
  gets $\rk(\gfrak^\theta)=\rk(\gfrak)$. This proves the proposition
  if $\gfrak$ is of type $B_n$, $C_n$, $E_7$, $E_8$, $F_4$, or $G_2$. 
  The other cases follow by inspection from 
  \cite[p.~514, Tables II and III]{b-Helga78} and Proposition \ref{PBZ}. For the
  convenience of the reader we collect the order $k$ of the diagram 
  automorphism $\nu$ for the remaining cases:
  
  \noindent{\bf \underline{$k=1$:}} $AIII$, $AIV$, ($DI$, case 1, $r$ even),
  ($DI$, case 2, $n$ odd), ($DI$, case 3, $n$ even), $DIII$, $EII$, $EIII$.
  
  \noindent{\bf \underline{$k=2$:}} $AI$, $AII$, ($DI$, case 1, $r$ odd),
  ($DI$, case 2, $n$ even), ($DI$, case 3, $n$ odd), $DII$, $EI$, $EIV$.
  $\blacksquare$

\vspace{.5cm}
\noindent{\bf 9.3 $Z(\Bc)$ is a polynomial ring in $\rk (\gfrak^\theta)$ variables.}
  Recall from (\ref{PBZ-free}) and Proposition \ref{rank=rank} that $\Ptheta$
  is generated over $\N$ by $m=\rk(\gfrak^\theta)$ linearly independent 
  elements $\nu_1,\dots,\nu_m\in \Ppiplus$. Let $d_{\nu_1},\dots,d_{\nu_m}$
  denote the corresponding central elements as chosen at the beginning
  of Section 8.4. Let $\cC[z_1,\dots,z_m]$ be the polynomial ring in
  $m$ generators.
  \setcounter{theorem}{2}
  \begin{theorem}
    The algebra homomorphism $\Phi:\cC[z_1,\dots,z_m]\rightarrow Z(\Bc)$
    defined by $\Phi(z_i)=d_{\nu_i}$ is an isomorphism.
  \end{theorem}
  {\bf Proof:} Recall from Theorem \ref{mainThm} that in the graded 
  algebra $\Gr$ the element $\odmu$ represented by $d_\mu$ lies in 
  $(\adr U)\otau$. Recall, moreover, 
  that for $\mu,\eta\in \Ptheta$ one has 
  $\odmu \,\odeta \in (\adr U)\overline{\tau(\mu+\eta)}$ 
  \cite[Prop.~4.12 (iii)]{a-JoLet2} and that $\Gr$ 
  is an integral domain \cite[Prop.~4.12 (iii)]{a-JoLet2}. Hence 
  $d_\mu d_\nu$ is a central element in $\Bc$ which consists of a 
  nonzero component in $(\adr U)\tau(\mu+\nu)$ and terms of lower 
  filter degree.
  By Theorem \ref{basis-prop} one gets 
  $d_\mu d_\eta=a_{\mu\eta}d_{\mu+\eta}$ for some constant
  $a_{\mu\eta}\in \cC$ up to terms of lower filter degree in $Z(\Bc)$.
  Hence up to terms of lower filter degree in $Z(\Bc)$ any element of
  $Z(\Bc)$ can be written as a polynomial in the generators 
  $d_{\nu_1},\dots,d_{\nu_m}$.  By induction this proves that $\Phi$ is
  surjective.
  
  As $\Gr$ is an integral domain the monomial
  $d_{\nu_1}^{n_1}\dots d_{\nu_m}^{n_m}$ represents a nonzero element
  in $(\adr U)\overline{\tau(2\sum n_i\nu_i)}$.
  Hence $\Phi$ is injective.
  $\blacksquare$

\section{Appendix: Commonly used  notation}

   \medskip
    \begin{tabbing}
	   \noindent
   \=Section 1.1:\quad\quad\quad\=\\
   \>$\C$\>    complex numbers\\
   \>$\Z$\>    integers\\
   \>$\N$\>    nonnegative integers \\
   \>$\gfrak$\> finite-dimensional complex semisimple Lie algebra\\
   \>$\hfrak$\> Cartan subalgebra of $\gfrak$\\
   \>$\Delta$\> root system of $\gfrak,\hfrak$ \\
   \>$\pi$\> set $\{\alpha_1,\dots, \alpha_n\}$ of simple roots for $\Delta$\\
   \>$(\cdot,\cdot)$\> Cartan inner product on $\hfrak^{\ast}$\\
   \>$W$\>   Weyl group of $\Delta$\\
   \>$w_0$\>   longest element in $W$ with respect to $\pi$\\
   \>$Q(\pi)$\>  root lattice of $\Delta$\\
   \>$P(\pi)$\>   weight lattice of $\Delta$\\
   \>$Q^+(\pi)$\> $\N\pi$\\
  \>$P^+(\pi)$\>  dominant integral weights associated to $\pi$\\
   \>$\omega_1.\dots,\omega_n$\> fundamental weights \\
   \>$\leq$\> standard partial ordering on $\hfrak^*$\\ 
  \>\>\\

   \>Section 1.2:\>\\
   \>$q$\> an indeterminate\\
   \>${\cal C}$\> ${\bf C}(q^{{1}\over{M}})$, $M$ a large integer\\
   \>$U_q(\gfrak)$\>  quantized enveloping algebra of $\gfrak$\\
   \>$t_i,x_i,y_i$\> generators of $U_q(\gfrak)$\\
   \>$U$\>$U_q(\gfrak)$\\
   \>$\kow$\>  coproduct of $U_q(\gfrak)$\\
   \>$\sigma$\>  antipode of $U_q(\gfrak)$\\
   \>$\epsilon$\>  counit of $U_q(\gfrak)$\\
   \>$\adr$\>  right adjoint action\\
   \>$U^+$\> subalgebra of $U$ generated by $x_1,\dots, x_n$\\
   \>$U^-$\> subalgebra of $U$ generated by $y_1,\dots, y_n$\\
   \>$G^+$\> subalgebra of $U$ generated by $x_1t_1^{-1},\dots, 
   x_nt_n^{-1}$\\
   \>$G^-$\> subalgebra of $U$ generated by $y_1t_1,\dots, y_nt_n$\\
   \>$T$\> group generated by $t_1^{\pm 1},\dots, t_n^{\pm 1}$\\
   \>$U^0$\> subalgebra of $U$ generated by $T$\\
   \>$\tau$\> group isomorphism from $Q(\pi)$ onto $T$\\
   \>$\check T$\> extension of $T$ isomorphic to $P(\pi)$\\
   \>$\check U$\> simply connected quantized enveloping algebra of 
   $\gfrak$\\
   \>$\check U^0$\> subalgebra of $\check U$ generated by 
   $\tau(\lambda),\lambda\in P(\pi)$\\
   \>$(yt)_I$\> $y_{i_1}t_{i_1}\dots y_{i_m}t_{i_m}$ where $I=(i_1,\dots, i_m)$\\
   \>${\rm wt}(I)$\>$\alpha_{i_1}+\dots+\alpha_{i_m}$ where $I=(i_1,\dots, i_m)$\\
   \>${\cal I}$\> set of multiindices such that $\{(yt)_I|\ I\in {\cal 
   I}\}$ is a basis for $G^-$.\\
   \>\>\\
   
   \>Section 1.3:\>\\
   \>$S_{\mu}$\> $\mu$ weight space of $S$\\
   \>$V(\mu)$\> simple left $U$-module of highest weight $\mu$\\
   \>$V(\mu)^*$\> dual space of $V(\mu)$ with natural right $U$-module structure\\
   \>\>\\
   
   \>Section 1.4:\>\\
   \>$\theta$\>involutive Lie algebra automorphism of $\gfrak$\\
   \>$\gfrak^{\theta}$\> $\{x\in \gfrak|\ \theta(x)=x\}$\\
   \>$\Theta$\> involution of $\hfrak^\ast$ induced by $\theta$\\
   \>$\pi_{\Theta}$\>$\{\alpha_i\in \pi|\Theta(\alpha_i)=\alpha_i\}$\\
   \>$p$\> permutation on $\{1,\dots, n\}$ satisfying (\ref{pdefn})\\
   \>$\pi^*$\>$\{\alpha_i\in \pi\setminus\pi_{\Theta}|\ i\leq p(i)\}$\\
   \>$\tilde\alpha$\>$(\alpha-\Theta(\alpha))/2$\\
   \>$\Sigma$\>  restricted root system $\{\tilde\alpha|\ \alpha\in 
   \Delta,\Theta(\alpha)\neq \alpha\}$\\
   \>$P(\Sigma)$\>  weight lattice of $\Sigma$\\
   \>$P^+(\Sigma)$\> set of dominant integral weights inside $P(\Sigma)$\\
   \>$\geq_r$\> standard partial ordering associated to restricted 
   root system\\
   \>\>\\
   
   \>Section 1.5:\>\\
   \>${\cal M}$\> algebra generated by $x_i,y_i,t_i^{\pm 
   1}$, for $\alpha_i\in \pi_{\Theta}$\\ 
   \>$T_{\Theta}$\> $\{\tau(\lambda)|\ \lambda\in Q(\pi)$ and 
   $\Theta(\lambda)=\lambda\}$\\
   \>$B$\>  algebra  generated by ${\cal M},T_{\Theta},B_i,$ 
     for $\alpha_i\in \pi\setminus\pi_{\Theta}$\\ 
     \>$B_i,\tilde\theta(y_i)$
     \>  for $\alpha_i\notin \pi_{\Theta}$, see (\ref{Bi-def})\\
   \>$B_i $\>  $y_it_i$ if $\alpha_i\in 
   \pi_{\Theta}$\\
   \>${\cal M}^+$\> ${\cal M}\cap U^+$\\
   \>$B_I$\>$B_{i_1}\dots B_{i_m}$ where $I=(i_1,\dots, i_m)$\\
   \>$T'_{\Theta}$\>$\{\tau(\lambda)|\ \lambda\in P(\pi)$ and 
   $\Theta(\lambda)=\lambda\}$\\
   \>$\check B$\> subalgebra of $\check U$ generated by $B$ 
   and $T'_{\Theta}$\\
   \>\>\\
   
   \>Section 1.6:\>\\
   \>$F_r(\check U)$\> $\{a\in \check U|\ {\rm dim}((\adr 
   U)a)<\infty\}$\\
   \>$Z(\check B)$\>  center of $\check B$\\
   \>\>\\
   
   \>Section 1.7:\>\\
   \>$M^B$\>  set of $B$-invariant elements in $M$\\
   \>spherical weight\>  weight in $2P^+(\Sigma)$\\
   \>\>\\
   
   \>Section 2.1:\>\\
   \>$h$\> height function on $\Q\pi$ defined by (\ref{h-def})\\
   \>${\cal F}$\> ${{1}\over{N}}\Z$-filtration defined by degree function 
   (\ref{filtdef})\\
   \>$\Gr$\> $\gr_{\cal F}(\check U)$\\
   \>$\Gr^0$\>$\gr_{\cal F}(\check U^0)$\\
   \>$\Gr$\> $\gr_{\cal F}(\check U)$\\
   \>$\Gr^-$\> $\gr_{\cal F}( U^-)$\\
   \>$\Gr^+$\> $\gr_{\cal F}(G^+)$\\
   \>$\bar{a}$\> graded image of $a$ in $\Gr$\\
   \>\>\\
   
   \>Section 2.2:\>\\
   \>$K(2\lambda)^-$\> subspace of $\Gr^-$ such that $(\adr 
   U^-)\overline{\tau(2\lambda)}=K(2\lambda)^-\otimes{\cal 
   C}\overline{\tau(2\lambda)}$\\
   \>$K(2\lambda)^+$\> subspace of $\Gr^+$ such that $(\adr 
     U^+)\overline{\tau(2\lambda)}=\Cc\overline{\tau(2\lambda)}\otimes K(2\lambda)^+$\\
   \>$F_r(\Gr)$\>$\{a\in \Gr|\ {\rm dim}((\adr U)a)<\infty\}$\\
   \>$g$\> isomorphism of $F_r(\check U)$ onto $F_r(\Gr)$ given by 
   (\ref{griso})\\
   \>\>\\
   
   \>Section 2.3:\>\\
   \>$h^{\Sigma}$\> height function associated to $\Sigma$ defined 
   by (\ref{restrictedheight})\\
   \>${\cal F}^{\theta}$\> ${{1}\over{N}}\Z$-filtration defined by degree 
   function (\ref{filtthetadefn})\\
   \>$\Gr^{\theta}$\>$\gr^{\theta}_{\cal F}(\check U)$\\
   \>$\bar{a}^{\theta}$\> graded image in $\Gr^{\theta}$ of $a\in 
   \check U$\\
   \>\>\\
   
   \>Section 2.4:\>\\
   \>$F_r(\Gr^{\theta})$\>$\{a\in \Gr^{\theta}|\ {\rm dim}((\adr U)a)<\infty\}$\\
   \>$g^{\theta}$\>isomorphism of $F_r(\check U)$ onto 
   $F_r(\Gr^{\theta})$ given in Proposition \ref{FGrtheta}\\
   \>\>\\
   
   \>Section 3.2:\>\\
   \>$W(\pi_{\Theta})$\> Weyl group for root system of $\pi_{\Theta}$\\
   \>$w_0'$\> longest element in $W(\pi_{\Theta})$\\
   \>$\mfrak$\> semisimple Lie subalgebra of $\gfrak$ with simple 
   roots $\pi_{\Theta}$\\
   \>\>\\
   
   \>Section 4.1:\>\\
   \>$T_{\geq}$\>${\rm span}\{\tau(\delta)|\ \delta\in P(\pi)$ and 
   $\tilde\delta\geq_r0\}$\\
   \>\>\\
   
   \>Section 4.2:\>\\
   \>$G({\le}\mu)$\>$U^-G^+\tau(\mu){\cal C}[t_1^{-2},\dots, 
   t_n^{-2}]$\\
   \>\>\\
   
   \>Section 5.1:\>\\
   \>$\check{\cal A}$\>$\{\tau(\tilde\alpha)|\ \alpha\in P(\pi)\}$\\
   \>$\Ttc$\>$\{\tau((\alpha+\Theta(\alpha))/2)\,|\, \alpha\in 
   P(\pi)\}$\\
   \>ad\> left adjoint action\\
   \>$N^+$\> subalgebra of $U^+$ generated by $(\ad {\cal M}^+)({\cal 
   C}[x_i\,|\, \alpha_i\notin \pi_{\Theta}])$\\
   \>\>\\
   
   \>Section 6.1:\>\\
   \>$N^+_+$\>intersection of $N^+$ with augmentation ideal of $U$\\
   \>$\check{\cal A}^+$\>$\check{\cal A}\setminus\{1\}$\\
   \>${\cal G}$\>$\check{\cal A}N^+_++\check{\cal A}^+$\\
   \>\>\\
   
   \>Section 6.2:\>\\
   \>$b$\> projection of $\check U\Ttc$ onto $\check B\Ttc$ defined by 
  (\ref{Iwoplus})\\
  \>\>\\
  
  \>Section 7:\>\\
  \>${\cal V}$\> $G^-{\cal M}^+T'_{\Theta}\cap F_r(\check 
  U)^{U^-}\cap \sum_{\lambda\in 2P^+(\Sigma)}\check U_{-\lambda}$\\
  \>\>\\
  
  \>Section 7.1:\>\\
  \>$\Ptheta$\>$\{\mu\in \Ppiplus|\ \Theta(\mu)=\mu+w_0\mu-w_0'\mu\}$\\
  \>\>\\
  
  \>Section 7.2:\>\\
  \>$Y^{\mu}$\> nonzero vector in ${\cal V}\cap (\adr 
  U)\tau(2\mu)$ of weight $2\tilde\mu$, \ $\mu\in \Ptheta$\\
  \>\>\\

  \>Section 7.3:\>\\
  \>$v_{\mu}^B$\> element in $((\adr U)Y^{\mu})^B$ defined by (\ref{vmuBdefn})\\
  \>$V^{\mu}_0$\> sum of simple spherical modules $V$ in
  $G({\le}2\mu)$ with
 $ V\not\cong V(2\tilde\mu)^*$\\
 \>\>\\
 
 \>Section 8.4:\>\\
 \>$d_{\mu}$\> element in $Z(\check B)\cap (v_{\mu}^B+(V_0^{\mu})^B)$, 
 for $\mu\in \Ptheta$\\
 \>\>\\
 
 \>Section 9.2:\>\\
 \>$\nu_1,\dots,\nu_m$\> linearly independent generators for $\Ptheta$
 (see (\ref{PBZ-free}))\\
 \>${\rm rank}(\Ptheta)$\>$m$ defined by (\ref{PBZ-free})\\ 
   
\end{tabbing}
   
\providecommand{\bysame}{\leavevmode\hbox to3em{\hrulefill}\thinspace}
\providecommand{\MR}{\relax\ifhmode\unskip\space\fi MR }
\providecommand{\MRhref}[2]{%
  \href{http://www.ams.org/mathscinet-getitem?mr=#1}{#2}
}
\providecommand{\href}[2]{#2}


\end{document}